# LOCAL PARTIAL-LIKELIHOOD ESTIMATION FOR LIFETIME DATA

By Jianqing Fan[1], Huazhen Lin and Yong Zhou[2]

*Chinese University of Hong Kong and Princeton University, Sichuan University and Chinese Academy of Science*

This paper considers a proportional hazards model, which allows one to examine the extent to which covariates interact nonlinearly with an exposure variable, for analysis of lifetime data. A local partial-likelihood technique is proposed to estimate nonlinear interactions. Asymptotic normality of the proposed estimator is established. The baseline hazard function, the bias and the variance of the local likelihood estimator are consistently estimated. In addition, a one-step local partial-likelihood estimator is presented to facilitate the computation of the proposed procedure and is demonstrated to be as efficient as the fully iterated local partial-likelihood estimator. Furthermore, a penalized local likelihood estimator is proposed to select important risk variables in the model. Numerical examples are used to illustrate the effectiveness of the proposed procedures.

**1. Introduction.** One of the most celebrated models for analyzing lifetime data is the Cox proportional hazards model, which explicitly postulates the covariate effects on the hazard risk via

$$\lambda(t) = \lambda_0(t) \exp\{g(\mathbf{Z})\},$$

where $\lambda_0(\cdot)$ is the baseline hazard risk and $g(\mathbf{Z})$ reflects the covariate effect. In parametric models it is commonly assumed that

$$g(\mathbf{Z}) = \boldsymbol{\beta}^T \mathbf{Z}$$

for some unknown parameters $\boldsymbol{\beta}$. See, for example, [1] and [20]. The log-linear model is a simple and mathematically convenient model that provides

Received December 2002; revised November 2004.

[1]Supported in part by NIH Grant R01 HL69720, NSF Grant DMS-03-55179 and RGC Grant CUHK 4262/01P of HKSAR.

[2]Supported in part by the Fund of National Natural Science (Grant 10471140) of China.

useful analysis for a covariate effect. However, in many biomedical studies, the covariate effects can be more complicated than the log-linear effect and new analytic challenges arise in assessing nonlinear effects. Beyond the traditional linear model, there are infinitely many possible nonlinear forms. Depending on the background of study, one often chooses a form that reasonably explains the objective of the study. For example, the effect of exposure variables and confounding factors on the hazard risk may vary with the level of an exposure variable, denoted by $W$. This leads one naturally to consider the model

$$\lambda(t) = \lambda_0(t) \exp\{\boldsymbol{\beta}(W(t))^T \mathbf{Z}(t) + g(W(t))\}. \tag{1.1}$$

Here $\boldsymbol{\beta}(\cdot)$ and $g(\cdot)$ are unknown coefficient functions, characterizing the extent to which the association varies with the level of the exposure variable $W$. Note that the term $g(W(t))$ can be incorporated into the covariates $\mathbf{Z}(t)$ by introducing a dummy variable with column one. We opt to not do so, because the local intercept for $g(\cdot)$ will cancel out in the local partial likelihood (2.3) below, leading to a different estimator rule for $g$. For ease of presentation, we drop the dependence of covariates on time $X_i$, with the understanding that the methods and proofs in this paper are applicable to time-dependent covariates.

When the variable $W$ is time, rather than a covariate variable, model (1.1) becomes a time-dependent coefficient Cox model, which has been studied by a number of authors, including Zucker and Karr [37], Murphy and Sen [31], Gamerman [21], Murphy [30], Marzec and Marzec [28], Martinussen, Scheike and Skovgaard [27], Cai and Sun [10], and Tian, Zucker and Wei [32]. In this case, unless the coefficient functions $\boldsymbol{\beta}(t)$ are independent of time $t$, the model is no longer a proportional hazards model. In contrast, model (1.1) is still a proportional hazards model. It allows one to examine the extent to which covariates $\mathbf{Z}$ interact nonlinearly with the exposure variable $W$. As will be explained later, although model (1.1) looks similar to the time-dependent coefficient Cox model, it is more involved when establishing asymptotic properties.

The varying-coefficient models arise from many different fields and have been studied in many different contexts. For cross-sectional type data, they have been studied as models to explore nonlinearity and assess nonlinear interactions by Cleveland, Grosse and Shyu [14], Hastie and Tibshirani [24], Carroll, Ruppert and Welsh [12], Fan and Zhang [19] and Cai, Fan and Li [8], among others. In time series, they are extensions of threshold autoregressive models and have been used to enhance the predictive power of linear autoregressive models. See, for example, [13] and [9]. The varying coefficient models have also been widely used to analyze longitudinal data. They allow one to examine the extent to which the association between independent and dependent variables varies over time. See, for example, [7, 25, 35, 36].

In this paper we propose techniques for estimating the coefficient functions $\boldsymbol{\beta}(\cdot)$ using local linear techniques [15]. The asymptotic bias and variance are obtained by establishing asymptotic normality. The variance is then estimated via a sandwich formula, which is shown to be consistent. To save computation of the local partial-likelihood estimator, a one-step procedure is proposed, which is shown to have the same asymptotic bias and variance as the local partial-likelihood estimator. Implementation of the proposed estimator depends on the choice of good initial estimators: estimates at the nearest grid points are recommended. The resulting procedure is demonstrated to be quite effective in our numerical implementation. In addition, the baseline hazard function $\lambda_0(\cdot)$ is estimated via a kernel method. The consistency property is demonstrated.

An objective of survival analysis is to identify the risk factors and their risk contributions. At the initial stage of a study, many covariates are collected to reduce possible modeling biases, and a large model is built, namely the dimensionality of $\mathbf{Z}$ in (1.1) is high. An important and challenging task is to efficiently select a subset of significant variables from model (1.1). Fan and Li [17] proposed a family of new variable selection methods based on a nonconcave penalized likelihood. Their methods are different from traditional ones in that they delete insignificant variables by estimating their coefficients as 0, and simultaneously select significant variables and estimate regression coefficients. Lasso, proposed by Tibshirani [33, 34], is a member of this family with an $L_1$ penalty. From their simulations, Fan and Li [17] showed that the penalized likelihood estimator with smoothly clipped absolute deviation (SCAD) penalty outperforms the best subset variable selection in terms of computational cost and stability in the terminology of Breiman [5]. In addition, they have proven that SCAD improves the lasso in terms of estimation biases. Furthermore, they have demonstrated that with a proper choice of regularization parameters and penalty functions (such as SCAD), the penalized likelihood estimator possesses an oracle property. Namely, the true regression coefficients that are zero are automatically estimated as zero and the remaining coefficients are estimated as well as if the correct submodel is known in advance. Hence, the SCAD and its siblings are ideal for variable selection, at least from a theoretical point of view. These nice properties encouraged us to extend the technique to the nonparametric model (1.1). It gives us a quick and effective method for eliminating unimportant variables.

The paper is organized as follows. Section 2 introduces the local partial-likelihood estimation and establishes the asymptotic normality. One-step estimation and estimation of the baseline hazard function are studied in Section 3. Section 4 deals with the issue of variable selection. Numerical examples are given in Section 5. Technical proofs are relegated to Appendix A.

**2. Partial-likelihood estimation.** Suppose that there is a random sample of size $n$ from an underlying population. Let $T_i$ denote the potential failure time, let $C_i$ denote the potential censoring time and let $X_i = \min(T_i, C_i)$ denote the observed time for the $i$th individual. Assume that $T_i$ and $C_i$ are independent given covariates $Z_i$ and $W_i$. Let $\Delta_i$ be an indicator which equals 1 if $X_i$ is a failure time and 0 otherwise. The covariates $\mathbf{Z}$ and $W$ are allowed to be time dependent. The observed data structure is

$$\{X_i, \Delta_i, \mathbf{Z}_i, W_i\} \qquad \text{for } i = 1, \dots, n,$$

where $\mathbf{Z}_i = (Z_{i1}, \dots, Z_{ip})^T$ and $W_i$ are two types of covariates, with $W$ being an exposure variable of interest.

When all the observations are independent, the partial likelihood for model (1.1) is

$$(2.1) \qquad L(\boldsymbol{\beta}(\cdot), g(\cdot)) = \prod_{i=1}^{n} \left\{ \frac{\exp\{\boldsymbol{\beta}(W_i)^T \mathbf{Z}_i + g(W_i)\}}{\sum_{j \in \mathcal{R}(X_i)} \exp\{\boldsymbol{\beta}(W_j)^T \mathbf{Z}_j + g(W_j)\}} \right\}^{\Delta_i},$$

where $R(t) = \{i : X_i \geq t\}$ denotes the set of the individuals at risk just prior to time $t$.

2.1. *Local partial likelihood.* If the unknown functions $\boldsymbol{\beta}(\cdot)$ and $g(\cdot)$ are parametrized, the parameters can be estimated by maximizing (2.1). For our nonparametric estimation, since the forms of the unknown functions are not available, we can only rely on their qualitative traits.

Assume that every component of $\boldsymbol{\beta}(\cdot)$ and $g(\cdot)$ is smooth so that it admits Taylor expansion: for each given $w_0$ and $w$ around $w_0$,

$$(2.2) \qquad \begin{aligned} \boldsymbol{\beta}(w) &\approx \boldsymbol{\beta}(w_0) + \boldsymbol{\beta}'(w_0)(w - w_0) \equiv \boldsymbol{\delta} + \boldsymbol{\eta}(w - w_0), \\ g(w) &\approx g(w_0) + g'(w_0)(w - w_0) \equiv \alpha + \gamma(w - w_0). \end{aligned}$$

Substituting this into (2.1), we obtain the logarithm of the local partial likelihood,

$$
\begin{aligned}
&\ell(\gamma, \boldsymbol{\delta}, \boldsymbol{\eta}) \\
&= n^{-1} \sum_{i=1}^{n} K_h(W_i - w_0) \Delta_i \\
(2.3) \quad &\times \Bigg\{ \boldsymbol{\delta}^T \mathbf{Z}_i + \boldsymbol{\eta}^T \mathbf{Z}_i(W_i - w_0) + \gamma(W_i - w_0) \\
&\qquad - \log\Bigg( \sum_{j \in \mathcal{R}(X_i)} \exp\{\boldsymbol{\delta}^T \mathbf{Z}_j + \boldsymbol{\eta}^T \mathbf{Z}_j(W_j - w_0) + \gamma(W_j - w_0)\} \\
&\qquad\qquad\qquad\qquad\qquad\qquad \times K_h(W_j - w_0) \Bigg) \Bigg\},
\end{aligned}
$$

where $K$ is a probability density called a kernel function, $h$ represents the size of the local neighborhood and $K_h(\cdot) = K(\cdot/h)/h$. The kernel weight is introduced to confirm that the local model (2.2) is only applied to the data around $w_0$. The local partial likelihood (2.3) can be derived from a profile likelihood point of view. The derivation is similar to those of Breslow [6] and Fan, Gijbels and King [16].

Let $\hat{\gamma}(w_0)$, $\hat{\boldsymbol{\delta}}(w_0)$ and $\hat{\eta}(w_0)$ be the maximizer of (2.3). Then $\hat{\boldsymbol{\beta}}(w_0) = \hat{\boldsymbol{\delta}}(w_0)$ is a local linear estimator for the coefficient function $\boldsymbol{\beta}(\cdot)$ at the point $w_0$. Similarly, an estimator of $g'(\cdot)$ at the point $w_0$ is simply the local slope $\hat{\gamma}(w_0)$, namely $\hat{g}'(w_0) = \hat{\gamma}(w_0)$. The curve $\hat{g}$ can be estimated by integration on the function $\hat{g}'(w_0)$. Following Hastie and Tibshirani [23], the integration can be approximated by using the trapezoidal rule.

We now express the local partial likelihood using the counting process notation. To this end, let $N_i(t) = I(T_i \leq t, \Delta_i = 1)$ and $Y_i(t) = I(X_i \geq t)$. Set

$$\boldsymbol{\xi} = (\boldsymbol{\delta}^T, \boldsymbol{\eta}^T, \gamma)^T \quad \text{and} \quad \mathbf{X}_i^* = (\mathbf{Z}_i^T, \mathbf{Z}_i^T(W_i - w_0), W_i - w_0)^T.$$

Then the local partial-likelihood function (2.3) can be expressed as

$$
\begin{aligned}
\ell_n(\boldsymbol{\xi}, \tau) = {} & n^{-1} \sum_{i=1}^n \int_0^\tau K_h(W_i - w_0) \boldsymbol{\xi}^T \mathbf{X}_i^* \, dN_i(u) \\
& - n^{-1} \sum_{i=1}^n \int_0^\tau K_h(W_i - w_0) \\
& \qquad \times \log \left\{ \sum_{j=1}^n Y_j(u) \exp(\boldsymbol{\xi}^T \mathbf{X}_j^*) K_h(W_j - w_0) \right\} dN_i(u)
\end{aligned}
$$
(2.4)

with $\tau = \infty$. To avoid the technicality of tail problems, only the data up to a finite time point $\tau$ are frequently used. Without ambiguity, we will let $\hat{\boldsymbol{\xi}}(w_0)$ be the maximizer of (2.4).

Note that the local partial likelihood in (2.4) is more complicated than that for the time-dependent coefficient Cox model. In particular, the kernel functions appear twice in the local partial likelihood (2.4), so as to use only local data. In contrast, for the time-dependent coefficient model, localizing in time once suffices. As a consequence, the technical proofs are more involved in the current setting.

The above method uses only one smoothing parameter to fit all the coefficient functions. When the coefficient functions admit different degree of smoothness [e.g., $g'(w)$ often admits a different degree of smoothness from other coefficient functions], one needs to use different bandwidths for different components. The two-step estimation method of Fan and Zhang [19] can be adapted here.

2.2. *Asymptotic normality.* We now establish the asymptotic normality of the local partial-likelihood estimator. As shown in Appendix A, the local partial-likelihood function $\ell_n(\boldsymbol{\xi}, \tau)$ is concave in $\boldsymbol{\xi}$ and its maximizer exists with probability tending to 1. Let $\mathbf{H}$ be a $(2p+1) \times (2p+1)$ diagonal matrix, with the first $p$ elements 1 and the remaining $p+1$ elements $h$, where $p$ is the number of elements in $\mathbf{Z}$. For any function $\boldsymbol{\xi}(w)$, $w \in J$, let $\|\boldsymbol{\xi}\|_J = \sup_{w \in J} |\boldsymbol{\xi}(w)|$, for a $p$-vector $\mathbf{a}$, let $|\mathbf{a}| = (\sum_{i=1}^p a_i^2)^{1/2}$ and $\|\mathbf{a}\| = \sup_i |a_i|$, and for a matrix $\mathbf{A}$, let $\|\mathbf{A}\| = \sup_{ij} |a_{ij}|$. Then we have the following consistency result.

THEOREM 1. *Under Conditions* A.1–A.8 *in Appendix* A, *we have*

$$\mathbf{H}\{\hat{\boldsymbol{\xi}}(w_0) - \boldsymbol{\xi}_0(w_0)\} \xrightarrow{P} 0,$$

*where* $\boldsymbol{\xi}_0(w_0) = (\boldsymbol{\beta}_0^T(w_0), \boldsymbol{\beta}_0'(w_0)^T, g_0'(w_0))^T$ *is the vector of the true parameter functions. If, in addition, Conditions* B.1–B.8 *hold, then we have the uniform consistency*

$$\|\mathbf{H}\{\hat{\boldsymbol{\xi}} - \boldsymbol{\xi}_0\}\|_{J_W} = \sup_{w \in J_W} |\mathbf{H}\{\hat{\boldsymbol{\xi}}(w) - \boldsymbol{\xi}_0(w)\}| \xrightarrow{P} 0,$$

*where* $J_W$ *is a compact subset of the support of the random variable* $W$.

To express explicitly the bias and variance of the estimator, we introduce some necessary notation. Let

$$\mu_i = \int x^i K(x)\, dx, \qquad \nu_i = \int x^i K^2(x)\, dx.$$

Denote

$$P(u, \mathbf{z}, w_0) = P(X \geq u | \mathbf{Z} = \mathbf{z}, W = w_0) \quad \text{and}$$

$$\rho(u, \mathbf{z}, w_0) = P(u, \mathbf{z}, w_0) \exp\{\boldsymbol{\beta}_0(w_0)^T \mathbf{z} + g_0(w_0)\}.$$

For $k = 0, 1, 2$, define

$$\mathbf{a}_k(u, w_0) = f(w_0) E\{\rho(u, \mathbf{Z}, w_0)\mathbf{Z}^{\otimes k} | W = w_0\},$$

where $f(\cdot)$ is the density of $W$ and $\mathbf{Z}^{\otimes k} = 1, \mathbf{Z}$ and $\mathbf{Z}\mathbf{Z}^T$ for $k = 0, 1$ and 2, respectively. Additionally set

$$\mathbf{a}_k = \mathbf{a}_k(w_0) = \int_0^\tau \mathbf{a}_k(u, w_0)\, d\Lambda_0(u).$$

We will drop the dependence of $\mathbf{a}_k(u, w_0)$ and $\mathbf{a}_k(w_0)$ on $w_0$ whenever there is no ambiguity. Finally, let

$$\boldsymbol{\Gamma} = \boldsymbol{\Gamma}(w_0) = \left\{\mathbf{a}_2 - \int_0^\tau \mathbf{a}_1(u)\mathbf{a}_1(u)^T \mathbf{a}_0^{-1}(u)\lambda_0(u)\, du\right\}^{-1}$$

and

$$\mathbf{Q} = \begin{pmatrix} (\mathbf{a}_2 - \mathbf{a}_1^T \mathbf{a}_0^{-1})^{-1} & -(\mathbf{a}_2 - \mathbf{a}_1^T \mathbf{a}_0^{-1})^{-1} \mathbf{a}_1 \mathbf{a}_0^{-1} \\ -\mathbf{a}_0^{-1} \mathbf{a}_1^T (\mathbf{a}_2 - \mathbf{a}_1 \mathbf{a}_1^T \mathbf{a}_0^{-1})^{-1} & (\mathbf{a}_0 - \mathbf{a}_1^T \mathbf{a}_2^{-1} \mathbf{a}_1)^{-1} \end{pmatrix},$$

where, in fact, $\mathbf{a}_0$ is a scale.

THEOREM 2 (Asymptotic normality). *Suppose that Conditions* A.1–A.8 *in Appendix* A *hold. Then*

$$\sqrt{nh}\{\mathbf{H}(\hat{\boldsymbol{\xi}}(w_0) - \boldsymbol{\xi}_0(w_0)) - \tfrac{1}{2}h^2 \mathbf{e}_p \boldsymbol{\xi}_0''(w_0)\mu_2\} \xrightarrow{\mathcal{L}} N(0, \boldsymbol{\Sigma}(\tau, w_0)),$$

*where* $\mathbf{e}_p$ *is a* $(2p+1)$-*order diagonal matrix, with the first* $p$ *elements* 1 *and the last* $p+1$ *diagonal elements* 0, $\boldsymbol{\xi}_0(w) = (\boldsymbol{\beta}_0^T(w), \boldsymbol{\beta}_0'(w)^T, g_0'(w_0))^T$ *and*

$$\boldsymbol{\Sigma}(\tau, w_0) = \begin{pmatrix} \boldsymbol{\Gamma}\nu_0 & \mathbf{0} \\ \mathbf{0}^T & \mathbf{Q}\mu_2^{-2}\nu_2 \end{pmatrix}.$$

The above theorem gives the joint asymptotic normality for the local partial-likelihood estimator. Its marginal distribution can easily be obtained as in the following corollary.

COROLLARY 1. *Under the conditions of Theorem* 2, *we have*

$$\sqrt{nh}\{\hat{\boldsymbol{\beta}}(w_0) - \boldsymbol{\beta}_0(w_0) - h^2 \boldsymbol{\beta}_0''(w_0)\mu_2/2\} \xrightarrow{\mathcal{L}} N(\mathbf{0}, \nu_0 \boldsymbol{\Gamma}),$$

$$\sqrt{nh^3}\{\hat{g}'(w_0) - g_0'(w_0)\} \xrightarrow{\mathcal{L}} N(0, (\mathbf{a}_0 - \mathbf{a}_1^T \mathbf{a}_2^{-1} \mathbf{a}_1)^{-1}\mu_2^{-2}\nu_2).$$

*Furthermore, they are asymptotically independent.*

As a consequence of Theorem 2, the theoretical optimal bandwidth can be obtained.

**3. Issues related to partial-likelihood estimation.** In this section we discussed a few issues that are related to the implementation of the partial-likelihood estimator.

3.1. *One-step local partial-likelihood estimator.* When estimating the whole functions $\boldsymbol{\beta}(\cdot)$ and $g(\cdot)$, we usually need to apply the local partial likelihood (2.4) at hundreds of points. Computing such an implicit estimator requires an iterative algorithm such as the Newton–Raphson method or Fisher's scoring method. Even worse, for certain given $w_0$, there does not exist a local partial-likelihood estimator due to the limited amount of data around $w_0$. These drawbacks make the local partial-likelihood estimator less appealing. Following Fan and Chen [15], we propose a one-step estimator as a viable alternative.

The local partial-likelihood estimator $\hat{\boldsymbol{\xi}}$ is found via solving the likelihood equation $\ell'_n(\boldsymbol{\xi}, \tau) = 0$, where $\ell'_n(\boldsymbol{\xi}, \tau) = \partial \ell_n(\boldsymbol{\xi}, \tau)/\partial \boldsymbol{\xi}$. To facilitate notation, from now on we drop the dependence of $\ell_n(\boldsymbol{\xi}, \tau)$ on $\tau$. For a given initial estimator $\hat{\boldsymbol{\xi}}_0$, by Taylor expansion we have

$$\ell'_n(\hat{\boldsymbol{\xi}}_0) + \ell''_n(\hat{\boldsymbol{\xi}}_0)(\hat{\boldsymbol{\xi}} - \hat{\boldsymbol{\xi}}_0) \approx 0.$$

Thus, the one-step estimator $\hat{\boldsymbol{\xi}}_{\mathrm{os}}$ is defined as

$$(3.1) \qquad \hat{\boldsymbol{\xi}}_{\mathrm{os}} = \hat{\boldsymbol{\xi}}_0 - \{\ell''_n(\hat{\boldsymbol{\xi}}_0)\}^{-1} \ell'_n(\hat{\boldsymbol{\xi}}_0).$$

A natural question arises: How good an initial estimator $\hat{\boldsymbol{\xi}}_0$ is needed for the one-step estimator to have the same performance as the maximum local partial-likelihood estimator. The following theorem gives an answer to this question.

THEOREM 3. *Under the conditions given in Theorem* 2, $\hat{\boldsymbol{\xi}}_{\mathrm{os}}$ *has the same asymptotic distribution as the maximum local partial-likelihood estimator* $\hat{\boldsymbol{\xi}}$, *provided that*

$$(3.2) \qquad \mathbf{H}(\hat{\boldsymbol{\xi}}_0 - \boldsymbol{\xi}_0) = O_p(h^2 + (nh)^{-1/2}).$$

Theorem 3 provides the conditions under which the one-step estimator performs as well as the local partial-likelihood estimator. However, it does not provide any guidance for choosing an initial estimator. Cai, Fan and Li [8] provided a useful strategy for the choice of initial estimators and their idea can be adapted to the current setting. The basic idea is first to compute the local partial-likelihood estimates at a few fixed points. Use these estimates as the initial values of their nearest grid points and obtain the one-step estimates at these grid points. For example, in our simulation studies we evaluate the functions at $n_{\mathrm{grid}} = 200$ grid points. We first compute the maximum local pseudo-partial-likelihood estimators at specific grid points $u_{20}, u_{60}, u_{100}, u_{140}$ and $u_{180}$, and then use them as the initial values for the one-step estimator at their nearest grid points. Use the newly computed one-step estimates (at points $u_{19}, u_{21}, u_{59}, u_{61}, \dots$) as the initial values of their nearest grid points to compute the one-step estimates and so on, until the one-step estimates at all grid points are computed. Hence, as long as the number of grid points is large enough, condition (3.2) holds.

3.2. *Estimation of baseline hazard function.* With estimators of $\boldsymbol{\beta}(\cdot)$ and $g(\cdot)$, we can estimate the baseline hazard function by using a kernel smoothing,

$$\hat{\lambda}_0(t) = \int W_b(t - x) \, d\hat{\Lambda}_0(x),$$

where $W$ is a given kernel function, $b$ is a given bandwidth and

$$\hat{\Lambda}_0(t) = \frac{1}{n} \sum_{i=1}^n \int_0^t \frac{dN_i(u)}{n^{-1} \sum_{j=1}^n Y_j(u) \exp(\hat{\boldsymbol{\beta}}^T(W_j) \mathbf{Z}_j(u) + \hat{g}(W_j))}.$$

Note that $\hat{\Lambda}_0(\cdot)$ is an estimate of the cumulative hazard function $\Lambda_0$.

THEOREM 4. *Under Condition* B *in Appendix* A, *we have*

$$\hat{\Lambda}_0(t) \longrightarrow \Lambda_0(t) \quad and \quad \hat{\lambda}_0(t) \longrightarrow \lambda_0(t)$$

*uniformly on* $(0, \tau]$ *in probability.*

3.3. *Estimation of biases and variances.* The biases of nonparametric estimates are generally hard to estimate, since they involve higher-order derivatives. However, their variances can be estimated quite reasonably. Thus, in construction of confidence intervals/bands, the bias components are frequently omitted; in particular, undersmoothing procedures have been used to make the biases negligible relative to their standard error. See, for example, [4, 22, 26]. Some people might argue that this is also the approach that parametric methods take—modeling biases are inevitable and they are simply ignored in the construction of the parametric confidence intervals.

The bias and covariance of these local estimators $H(\hat{\boldsymbol{\xi}}(w_0) - \boldsymbol{\xi}_0(w_0))$ can be estimated by

$$\hat{\mathbf{A}}_n^{-1}(\tau, w_0) \hat{\mathbf{B}}_n(\tau, w_0) \quad \text{and} \quad (nh)^{-1} \hat{\mathbf{A}}_n^{-1}(\tau, w_0) \hat{\mathbf{\Pi}}_n(\tau, w_0) \hat{\mathbf{A}}_n^{-1}(\tau, w_0),$$

where

$$\hat{\mathbf{A}}_n(\tau, w_0) = \frac{1}{n} \sum_{i=1}^n \int_0^\tau K_h(W_i - w_0)$$

$$\times \frac{\hat{S}_{n2}(u, w_0) \hat{S}_{n0}(u, w_0) - \hat{S}_{n1}(u, w_0) \hat{S}_{n1}(u, w_0)^T}{(\hat{S}_{n0}(u, w_0))^2} \, dN_i(u),$$

$$\hat{\mathbf{B}}_n(\tau, w_0) = \frac{1}{n} \sum_{i=1}^n \int_0^\tau K_h(W_i - w_0) \left( \mathbf{U}_i^*(u) - \frac{\hat{S}_{n1}(u, w_0)}{\hat{S}_{n0}(u, w_0)} \right) Y_i(u) \hat{\lambda}_i(u) \, du,$$

$$\hat{\mathbf{\Pi}}_n(\tau, w_0) = \frac{h}{n} \sum_{i=1}^n \int_0^\tau K_h^2(W_i - w_0) \left( \mathbf{U}_i^*(u) - \frac{\hat{S}_{n1}(u, w_0)}{\hat{S}_{n0}(u, w_0)} \right)^{\otimes 2} Y_i(u) \hat{\lambda}_i(u) \, du,$$

with $\hat{\lambda}_i(u) = \exp(\hat{\boldsymbol{\beta}}(W_i)^T \mathbf{Z}_i(u) + \hat{g}(W_i)) \hat{\lambda}_0(u)$, $\mathbf{U}_i^* = \mathbf{H}^{-1} \mathbf{X}_i^*$ and

$$\hat{S}_{nk}(u, w_0) = \sum_{i=1}^n K_h(W_i - w_0) Y_i(u) \exp(\hat{\boldsymbol{\xi}}_0^T(w_0) \mathbf{X}_i^*(u)) (\mathbf{U}_i^*(u))^{\otimes k},$$

$$k = 0, 1, 2.$$

THEOREM 5. *Under the conditions of Theorem 4, we have*

$$h^{-2}\hat{\mathbf{A}}_n^{-1}(\tau, w_0)\hat{\mathbf{B}}_n(\tau, w_0) \longrightarrow \mathbf{e}_p\boldsymbol{\xi}''(w_0)\mu_2/2,$$

$$\hat{\mathbf{A}}_n^{-1}(\tau, w_0)\hat{\boldsymbol{\Pi}}_n(\tau, w_0)\hat{\mathbf{A}}_n^{-1}(\tau, w_0) \longrightarrow \boldsymbol{\Sigma}(\tau, w_0)$$

*in probability.*

In fact, by using the martingale properties, we can construct different estimators of $\hat{\mathbf{B}}_n(\tau, w_0)$ and $\hat{\boldsymbol{\Pi}}(\tau, w_0)$ without estimating the baseline hazard function $\lambda_0(\cdot)$. That is,

$$\widetilde{\mathbf{B}}_n(\tau, w_0) = \frac{1}{n}\sum_{i=1}^n \int_0^\tau K_h(W_i - w_0)\left(\mathbf{U}_i^*(u) - \frac{\hat{S}_{n1}(u, w_0)}{\hat{S}_{n0}(u, w_0)}\right) dN_i(u),$$

$$\widetilde{\boldsymbol{\Pi}}_n(\tau, w_0) = \frac{h}{n}\sum_{i=1}^n \int_0^\tau K_h^2(W_i - w_0)\left(\mathbf{U}_i^*(u) - \frac{\hat{S}_{n1}(u, w_0)}{\hat{S}_{n0}(u, w_0)}\right)^{\otimes 2} dN_i(u).$$

The results of Theorem 5 still hold when the quantities $\hat{B}_n(\tau, w_0)$ and $\hat{\boldsymbol{\Pi}}_n(\tau, w_0)$ are replaced by $\widetilde{B}_n(\tau, w_0)$ and $\widetilde{\boldsymbol{\Pi}}_n(\tau, w_0)$, respectively.

One can also use the bootstrap method as in [32] to obtain an estimated variance for our estimators. In fact, the method is particularly useful for estimating the sampling variability of $\hat{g}(w)$, since its analytic form is hard to derive.

## 4. Variable selection via nonconcave penalized likelihood.

4.1. *Local penalized likelihood.* For the nonparametric model (1.1), it is not easy to give a variable selection procedure without going to detailed inferences on each coefficient function. Motivated by the work of Fan and Li [17, 18], we apply their procedure locally around each grid point $w_0$. This results in the penalized log partial-likelihood function

$$(4.1) \qquad Q(\boldsymbol{\xi}) = \ell_n(\boldsymbol{\xi}, \tau) - \sum_{j=1}^{2p+1} p_\varrho(|\xi_j|),$$

where $p_\varrho(\cdot)$ is a penalty function. The penalized local partial-likelihood estimate of $\boldsymbol{\xi}$ is to maximize (4.1). With a proper choice of $\varrho$ and a penalty function, many estimated coefficients will be zero and hence their corresponding variables do not appear in the model at the point $w_0$. This achieves the objective of variable selection and results in a simple and implementable method to begin with.

A good penalty function should result in an estimator with the following three properties: unbiasedness for large coefficients to attenuate biases, sparsity (many small coefficients are estimated as zero) to reduce model complexity and continuity to avoid unnecessary variation in model prediction. Necessary conditions for unbiasedness, sparsity and continuity have been derived by Antoniadis and Fan [3] and Fan and Li [17]. A simple penalty function that satisfies all the three mathematical requirements is the smoothly clipped absolute deviation (SCAD) penalty, defined by

$$p'_\varrho(\theta) = \varrho \left\{ I(\theta \le \varrho) + \frac{(a\varrho - \theta)_+}{(a-1)\varrho} I(\theta > \varrho) \right\}$$

(4.2)

$$\text{for some } a > 2 \text{ and } \theta > 0.$$

Fan and Li [17] suggested using $a = 3.7$ from a Bayesian point of view and this value will be used in our numerical implementation.

There are two issues related to the practical implementation of the procedure. First, to facilitate the implementation we use only one regularization parameter for all variables which can have very different scales. Thus, we need to standardize variables before using (4.1). Since each variable in (4.1) is used locally around a given point $w_0$, its sample mean and standard deviation should be defined locally. For example, the variable $Z_1$ at the point $w_0$ can be standardized by

$$\text{ave}(Z_1|w_0) = \frac{1}{N} \sum_{i=1}^{n} K_h(W_i - w_0) Z_{1i}$$

and

$$\text{var}(Z_1|w_0) = \frac{1}{N} \left( \sum_{i=1}^{n} K_h(W_i - w_0) Z_{1i}^2 - \text{ave}(Z_1|w_0)^2 N \right),$$

where $N = \sum_{i=1}^{n} K_h(W_i - w_0)$. The second issue is that the number of variables as a function of $w_0$, if not constant, will be discontinuous. This will lead to discontinuous estimates of coefficient functions. This may not be bad in terms of overall prediction error, but does not produce parsimonious and appealing models. To avoid this, we use a simple voting rule: if a coefficient function is estimated as zero over a certain percentage of grid points, delete its corresponding variable; otherwise keep the variable. In our implementation, we use the majority voting rule, namely, the thresholding percentage is taken as 50%.

### 4.2. *Oracle property.*

We now establish an oracle property of the penalized local partial-likelihood estimator. We assume without loss of generality that the first $s$ variables of $\mathbf{Z}$ are significant and the last $p - s$ variables are not significant. To state our main result more explicitly, we need the following notation.

Recall that $\boldsymbol{\xi} = (\boldsymbol{\delta}^T, \boldsymbol{\eta}^T, \gamma)^T$. We divide $\boldsymbol{\delta}$ into $(\boldsymbol{\delta}_1^T, \boldsymbol{\delta}_2^T)^T$, where $\boldsymbol{\delta}_1$ and $\boldsymbol{\delta}_2$ are $s \times 1$ and $(p-s) \times 1$ vectors, representing, respectively, the vanishing and nonvanishing coefficients. Corresponding to the partition of $\boldsymbol{\delta}$, we divide $\boldsymbol{\eta}$ into $(\boldsymbol{\eta}_1^T, \boldsymbol{\eta}_2^T)^T$. Write

$$\boldsymbol{\xi}_1 = (\boldsymbol{\delta}_1^T, \boldsymbol{\eta}_1^T, \gamma)^T = (\xi_{1,1}, \ldots, \xi_{1,2s}, \xi_{1,2s+1})^T$$

and $\boldsymbol{\xi}_2 = (\boldsymbol{\delta}_2^T, \boldsymbol{\eta}_2^T)^T$. Let $\boldsymbol{\xi}_{10} = (\xi_{1,1,0}, \ldots, \xi_{1,2s+1,0})^T$, and $\boldsymbol{\xi}_{20}$ and $\boldsymbol{\xi}_0$ be, respectively, the true values of $\boldsymbol{\xi}_1$, $\boldsymbol{\xi}_2$ and $\boldsymbol{\xi}$. For example, $\xi_{1,j,0} = \beta_{j0}(w_0)$ for $j = 1, \ldots, s$, $\xi_{1,j,0} = \beta'_{j0}(w_0)$ for $j = s+1, \ldots, 2s$ and $\xi_{1,2s+1,0} = g'_0(w_0)$. Without loss of generality, assume that $\boldsymbol{\xi}_{20} = 0$. Set

$$a_n(w_0) = \max\{p'_\varrho(|\xi_{1,j,0}|) : \xi_{1,j,0} \neq 0\},$$

$$b_n(w_0) = \max\{p''_\varrho(|\xi_{1,j,0}|) : \xi_{1,j,0} \neq 0\}.$$

Let $\boldsymbol{\Pi}_1$ and $\mathbf{A}_1$ be, respectively, the submatrices of $\boldsymbol{\Pi}(\tau, w_0)$ and $\mathbf{A}(\tau, w_0)$ in (A.10) and (A.16) in Appendix A that correspond to the rows in $\boldsymbol{\xi}_1$. Corresponding to the partition of $\boldsymbol{\delta}$, let $\boldsymbol{\Gamma}^{-1} = (\boldsymbol{\Gamma}_{-1}^T, \boldsymbol{\Gamma}_{-2}^T)^T$ with $\boldsymbol{\Gamma}_{-1}$ and $\boldsymbol{\Gamma}_{-2}$ being $s \times p$ and $(p-s) \times p$ matrices, respectively.

The following theorem shows how the rates of convergence for the penalized local partial-likelihood estimates depend on the regularization parameter.

THEOREM 6. *Suppose that Conditions* A.1–A.8 *in the Appendix* A *hold. If* $b_n(w_0) \to 0$, *then there exists a local maximizer* $\hat{\boldsymbol{\xi}}_p$ *of* $Q(\boldsymbol{\xi})$ *such that* $\|\hat{\boldsymbol{\xi}}_p - \boldsymbol{\xi}_0\| = O_p(h^2 + (nh)^{-1/2} + a_n(w_0))$.

It is clear from Theorem 6 that by choosing a proper $\varrho$, such that $a_n(w_0) = O((nh)^{-1/2} + h^2)$, there exists a $(nh)^{-1/2} + h^2$ consistent penalized local partial-likelihood estimator. Now we show that this estimator must possess an oracle property.

THEOREM 7. *Assume that the penalty function* $p_\varrho(\theta)$ *satisfies*

$$(4.3) \qquad \liminf_{n \to \infty} \liminf_{\theta \to 0+} p'_\varrho(\theta)/\varrho > 0.$$

*Let* $\varrho \to 0$, $\{(nh)^{-1/2} + h^2\}/\varrho \to 0$ *and* $a_n(w_0) = O((nh)^{-1/2} + h^2)$. *Under the conditions of Theorem* 6, *the consistent local maximizer* $\hat{\boldsymbol{\xi}}_p = (\hat{\boldsymbol{\xi}}_{1p}^T, \hat{\boldsymbol{\xi}}_{2p}^T)$ *in Theorem* 6 *satisfies the following statements with probability tending to* 1:

(a) (Sparsity) *We have* $\hat{\boldsymbol{\xi}}_{2p} = \mathbf{0}$.
(b) (Asymptotic normality) *We have*

$$(4.4) \qquad \sqrt{nh}\mathbf{B}_1 \bigg\{ H_1(\hat{\boldsymbol{\xi}}_{1p} - \boldsymbol{\xi}_{10}) \\ - \mathbf{B}_1^{-1}\bigg[\mathbf{H}_1^{-1}\mathbf{b} + h^2\mu_2 \begin{pmatrix} \boldsymbol{\Gamma}_{-1}\boldsymbol{\beta}''_0(w_0) \\ \mathbf{0}_{s+1} \end{pmatrix} \bigg] \bigg\} \longrightarrow N(\mathbf{0}, \boldsymbol{\Pi}_1(\tau, w_0)),$$

*where* $\mathbf{b} = (p'_\varrho(|\xi_{1,1,0}|)\operatorname{sgn}(\xi_{1,1,0}),\ldots,p'_\varrho(|\xi_{1,2s+1,0}|)\operatorname{sgn}(\xi_{1,2s+1,0}))^T$, $\mathbf{B}_1 = \mathbf{A}_1 - \mathbf{H}_1^{-1}\boldsymbol{\Sigma}_1\mathbf{H}_1^{-1}$, $\boldsymbol{\Sigma}_1 = \operatorname{diag}\{p''_\varrho(|\xi_{1,1,0}|),\ldots,p''_\varrho(|\xi_{1,2s+1,0}|)\}$, $\boldsymbol{\beta}_0(w_0) = (\beta_{10}(w_0),\beta_{20}(w_0),\ldots,\beta_{s0}(w_0),\mathbf{0})^{\mathbf{T}}$, *and* $\mathbf{H}_1$ *is a* $(2s+1)\times(2s+1)$ *diagonal matrix with first* $s$ *elements* 1 *and the last* $s+1$ *elements* $h$.

We now explain that the penalized local-likelihood estimators possess an oracle property when penalty functions are properly chosen. Suppose that there is an oracle who knows $\boldsymbol{\xi}_{2p} = 0$. She then uses this knowledge to estimate $\hat{\boldsymbol{\xi}}_{1p}$, resulting in an oracle estimator. From Theorem 2, the asymptotic covariance matrix of this oracle estimator is $\frac{1}{nh}\mathbf{A}_1^{-1}\boldsymbol{\Pi}_1(\tau,w_0)\mathbf{A}_1^{-1}$. For penalty functions such as SCAD, since $\varrho \to 0$, for sufficiently large $n$,

$$a_n(w_0) = 0 \quad \text{and} \quad b_n(w_0) = 0 \qquad \text{so} \quad \mathbf{b} = 0 \quad \text{and} \quad \boldsymbol{\Sigma}_1 = 0.$$

Thus, Theorems 6 and 7 yield that $\hat{\boldsymbol{\xi}}_{2p} = 0$ and $H_1(\hat{\boldsymbol{\xi}}_{1p} - \boldsymbol{\xi}_{10})$ is asymptotically normal with covariance matrix $\frac{1}{nh}\mathbf{A}_1^{-1}\boldsymbol{\Pi}_1(\tau,w_0)\mathbf{A}_1^{-1}$, which is the same as the asymptotic variance of the oracle estimator (see Theorem 2). Furthermore, it can easily be seen that both estimators share the same asymptotic bias. Thus, the penalized likelihood estimators perform as well as the oracle estimator when the penalty functions are constant at the tails. In other words, when the true parameters have some zero components, they are estimated as 0 with probability tending to 1 and the nonzero components are estimated as well as the case where the correct submodel is known.

## 5. Numerical examples.

5.1. *Simulations.* In this section we first compare the performance of the one-step and local partial-likelihood estimators. The performance of estimator $\hat{\boldsymbol{\beta}}(\cdot)$ is assessed via the weighted mean square error (WMSE),

$$(5.1) \qquad \text{WMSE} = \frac{1}{n_{\text{grid}}}\sum_{j=1}^{p}\sum_{k=1}^{n_{\text{grid}}}a_j[\hat{\beta}_j(w_k) - \beta_j(w_k)]^2,$$

or the unweighted mean square error (UMSE) with all $a_j = 1$, where $\{w_k, k = 1,\ldots,n_{\text{grid}}\}$ are the grid points at which the functions $\boldsymbol{\beta}(\cdot)$ are estimated. In the following examples, the Gaussian kernel will be used, $n_{\text{grid}} = 200$ and, for WMSE, $a_j$ is reciprocal to the sample variance of $\{\beta_j(w_k)\}$.

EXAMPLE 1. We first consider the varying-coefficient model $\lambda(t) = 4t^3 \times \exp\{b(Z_1(t), Z_2, W)\}$ with

$$b(Z_1, Z_2, W) = 0.5W(1.5 - W)Z_1 + \sin(2W)Z_2$$
$$+ 0.5\{\exp(W - 1.5) - \exp(-1.5)\},$$

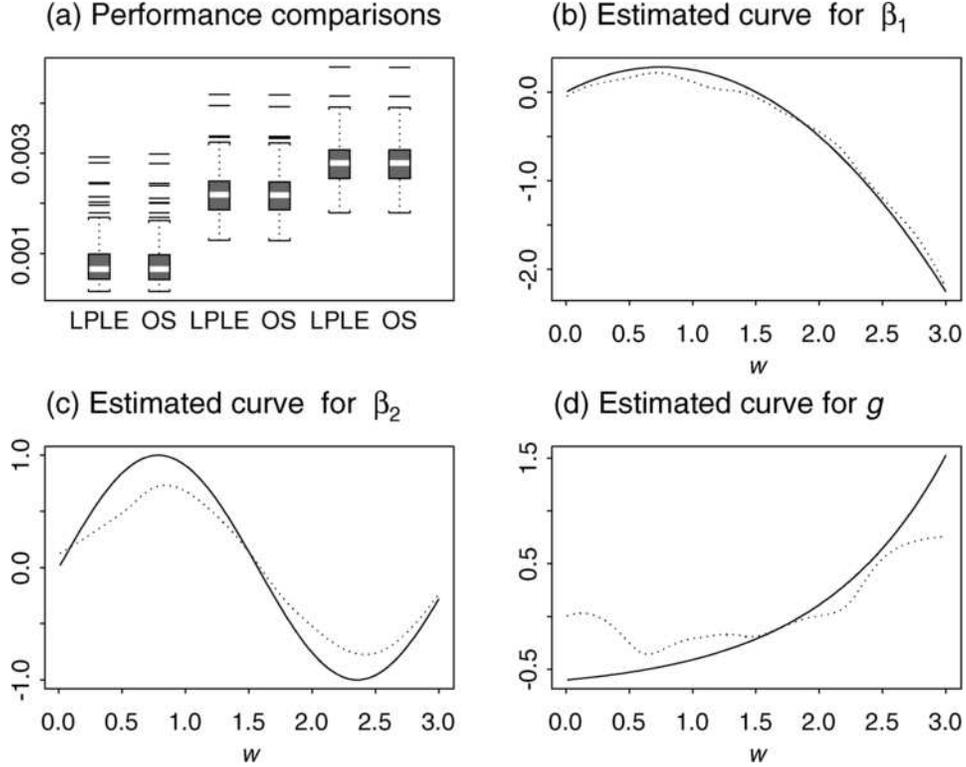

Fig. 1. *Simulation results for Example* 1. (a) *Boxplots for the distribution for the WMSE over the* 200 *replications, using the three bandwidths* $h = 0.2, 0.5, 1$ *(from left to right).* (b), (c) *and* (d) *Typical estimates of* $\beta_1(\cdot)$, $\beta_2(\cdot)$ *and* $g(\cdot)$ *with bandwidth* $h = 0.2$ *(solid line, true function; dashed line, one-step LPLE, i.e., OS).*

where $W$ is a random variable uniformly distributed on $[0, 3]$, the covariate $Z_1(t)$ is time-dependent, defined as $Z_1(t) = Z_1/4I(t \leq 1) + Z_1 I(t > 1)$, and $Z_1$ and $Z_2$ are jointly normal with correlation 0.5, each with mean 0 and standard deviation 5. The censoring random variable $C$ given $(Z_1, Z_2, W)$ is distributed uniformly on $[0, a(Z_1, Z_2, W)]$, where

$$a(Z_1, Z_2, W) = c_1 I(b(Z_1, Z_2, W) > b_0) + c_2 I(b(Z_1, Z_2, W) \leq b_0),$$

with $b_0$ being the mean function of $b(Z_1, Z_2, W)$. The constants $c_1 = 0.8$ and $c_2 = 20$ are chosen so that about 30–40% of data are censored in each region of the function $a(\cdot)$.

We have conducted 200 simulations with sample size 300. Figure 1(a) depicts the distribution for the WMSE over the 200 replications, using the three bandwidths $h = 0.2, 0.5, 1$. The initial value is chosen at grid points $w_{20}, w_{60}, w_{100}, w_{140}$ and $u_{180}$ by the local partial-likelihood estimator just mentioned in Section 3.1. It is evident that the performances of the one-step local partial-likelihood estimator (one-step LPLE) and local partial-likelihood estimator (LPLE) are comparable for a wide range of bandwidths. Figure 1(b)–(d) presents estimates of the coefficient functions from a typical sample (attaining the median WMSE performance) with $h = 0.2$.

We now test the accuracy of our standard error formula given in Section 3.3. The standard deviations, denoted by SD in Table 1, of 200 estimated $\hat{\beta}_1(w_0)$, $\hat{\beta}_2(w_0)$ and $\hat{g}'(w_0)$, based on 200 simulations, can be regarded as the true standard errors. The average and the standard deviation of 200 estimated standard errors, denoted by $\mathrm{SE}_{\mathrm{ave}}$ and $\mathrm{SE}_{\mathrm{std}}$, summarize the overall performance of the standard error formula. Table 1 presents the results at the points $w = 0.3, 0.75, 1.5, 2.25$ and $2.7$, which correspond to the 10th, 25th, 50th, 75th and 90th percentiles of the distribution of $W$. The performance of the standard error formula is quite satisfactory.

EXAMPLE 2.   In the following examples, we evaluate the performance of the proposed variable selection method. Samples of size 300 were simulated from the hazard regression model

$$\lambda(t) = \exp\left(\sum_{j=1}^{4} Z_j \beta_j(w) + g(w)\right),$$

where $\beta_1(w) = 3(w-2)^2$, $\beta_2(w) = 4\cos(\frac{(w-1.5)\pi}{5})$ and $\beta_3(w) = \beta_4(w) = g(w) = 0$. The covariates $Z_1, Z_2, Z_3$ and $Z_4$ are jointly normal, all with mean 0 and variance 2, and pairwise correlation 0.6. They are independent of $W$, which is uniformly distributed on $[0, 3]$. The censoring time follows the uniform distribution on $[0, 7]$ so that about 30–40% of the data were censored. The kernel function is Gaussian.

The performance of the proposed variable selection technique is compared with that of the maximum local partial-likelihood estimator from the full model and from the oracle estimator, which is based on the model with only covariates $Z_1$ and $Z_2$. Figure 2 depicts the distribution for the UMSE over the 200 replications, using bandwidths $h = 0.3$ and $\lambda = 0.3$. It is evident that the proposed variable selection procedure outperforms the maximum local partial-likelihood estimator and performs comparably with the oracle estimator.

TABLE 1
*True and estimated standard errors using bandwidth = 0.2 for Example 1*

| $w_0$ | $\hat{\beta}_1(w_0)$ | | | $\hat{\beta}_2(w_0)$ | | | $\hat{g}'(w_0)$ | | |
|---|---|---|---|---|---|---|---|---|---|
| | SD | $\mathrm{SE}_{\mathrm{ave}}$ | $(\mathrm{SE}_{\mathrm{std}})$ | SD | $\mathrm{SE}_{\mathrm{ave}}$ | $(\mathrm{SE}_{\mathrm{std}})$ | SD | $\mathrm{SE}_{\mathrm{ave}}$ | $(\mathrm{SE}_{\mathrm{std}})$ |
| 0.30 | 0.0606 | 0.0573 | (0.0098) | 0.0655 | 0.0479 | (0.0111) | 0.3831 | 0.3735 | (0.0492) |
| 0.75 | 0.0458 | 0.0479 | (0.0076) | 0.0579 | 0.0337 | (0.0079) | 0.2779 | 0.2967 | (0.0354) |
| 1.50 | 0.0340 | 0.0414 | (0.0058) | 0.0473 | 0.0236 | (0.0043) | 0.1910 | 0.2457 | (0.0258) |
| 2.25 | 0.0303 | 0.0343 | (0.0046) | 0.0282 | 0.0197 | (0.0018) | 0.1873 | 0.1602 | (0.0228) |
| 2.70 | 0.0429 | 0.0385 | (0.0053) | 0.0321 | 0.0222 | (0.0027) | 0.2491 | 0.1474 | (0.0178) |

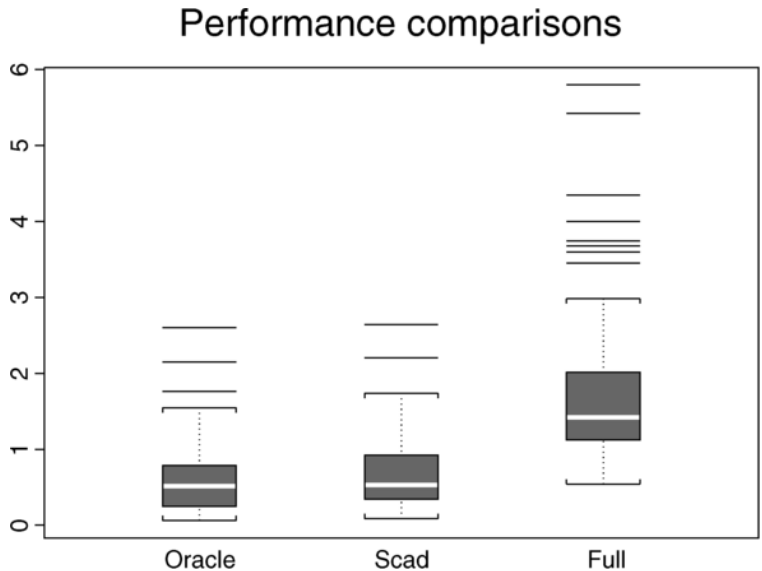

Fig. 2. *Boxplot for the distribution of the UMSE over the* 200 *replications, using bandwidths* $h = 0.3$ *and* $\lambda = 0.3$.

Using the majority voting (50%) rule, the variables $Z_3$, $Z_4$ and $g(W)$ were simultaneously deleted 98.5% of the time among 200 simulations, and using a 60% thresholding level, the variables $Z_3, Z_4$ and $g(w)$ were simultaneously deleted 92% of the time. Hence, only variables $Z_1$ and $Z_2$ remain. Their estimated coefficients are depicted in Figure 3 for a typical sample.

5.2. *Data analysis.* The proposed approaches are now applied to the nursing home data set analyzed by Morris, Norton and Zhou [29], where a full description of this data set is given. The data are from an experiment sponsored by the National Center for Health Services Research during 1980–1982 that involved 36 for-profit nursing homes in San Diego, California, with a sample of size 1601.

The study was designed to evaluate the effects of different financial incentives on, among other things, the duration of stay. This motivated Morris, Norton and Zhou [29] to take days $T$ in the nursing home as the response

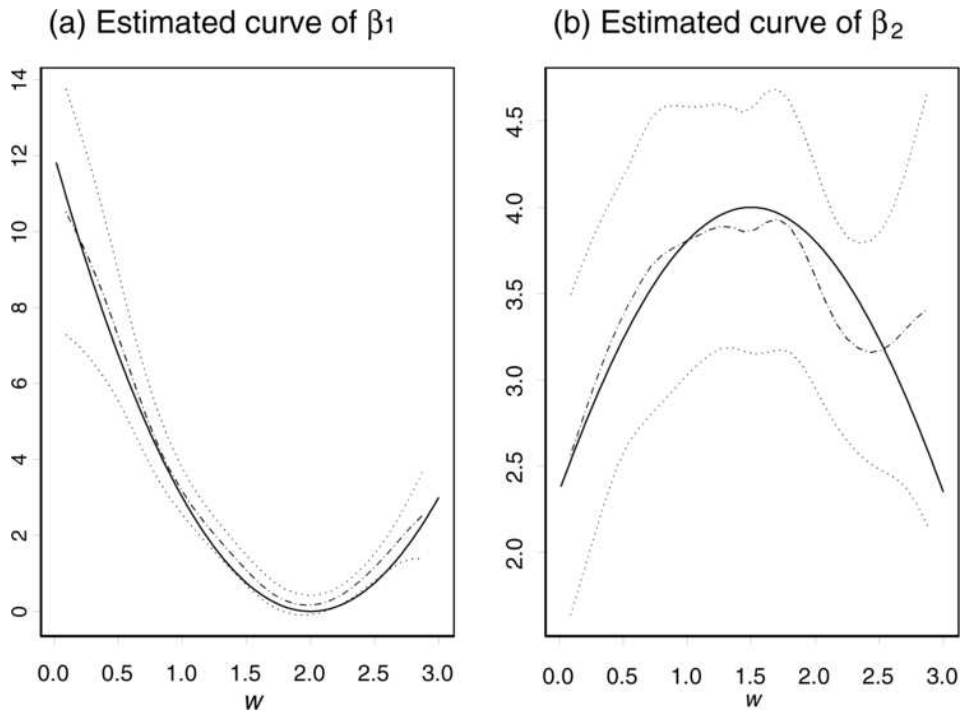

**(a) Estimated curve of $\beta_1$**

**(b) Estimated curve of $\beta_2$**

Fig. 3.  *The estimated coefficient functions (dashed lines) using the local partial-likelihood approach with bandwidth $h = 0.3$ after deleting $Z_3, Z_4$ and $g(w)$, as well as true lines (solid lines) and their 95% confidence bands (dotted lines) for Example 2.*

variable. They used the model

$$\lambda(t, x) = \lambda_0(t) \exp\left( \sum_{j=1}^{7} x_j \beta_j \right),$$

where $x_1$ is a treatment indicator, being 1 if treated at a nursing home and 0 otherwise; $x_2$ is a gender variable (1 for males and 0 for females); $x_3$ is a marital status indicator (1 if married and 0 otherwise); $x_4, x_5, x_6$ are three binary health status indicators that correspond to the best health to the worst health; $x_7$ is age, which ranges from 65 to 104. Morris, Norton and Zhou [29] fitted the Cox model with three parametric and one nonparametric baseline hazard model to this data set. Their model does not include any possible interactions between age and other variables. To explore possible interaction, Fan and Li [17] added interaction terms such as $x_7 x_1, x_7 x_2, \ldots$ in the initial model. With our newly developed technique, we can fit the

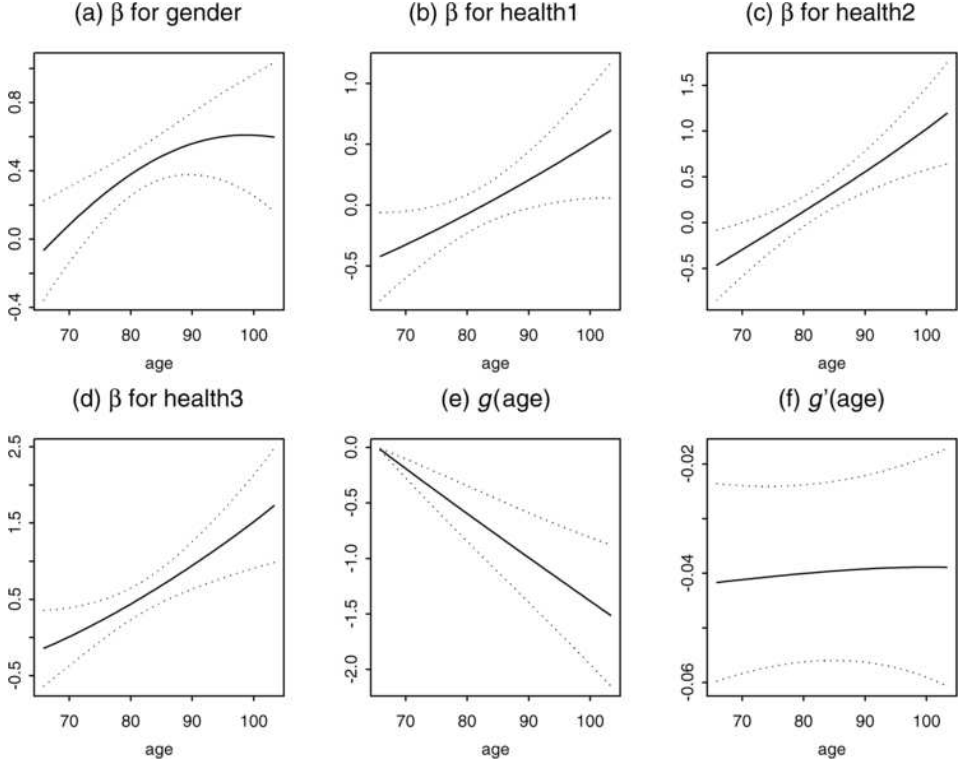

Fig. 4. *The estimated coefficient functions (solid lines) via a local partial-likelihood approach with bandwidth $h = 15$ and their 95% confidence limits (dotted lines) for the nursing home data without the treatment and marital covariates.*

more general model

$$\lambda(t, x) = \lambda_0(t) \exp\left( \sum_{j=1}^{6} \beta_j(x_7) x_j + g(x_7) \right).$$

This permits us to examine how different age groups interact with covariates such as treatment, gender and marital status. In fact, as age increases, elderly people would expect to stay at nursing homes longer. Therefore, it is natural to introduce the term $g(x_7)$, the varying intercept.

The local partial-likelihood method was applied to the data set with bandwidth $h = 15$, which was chosen by $K$-fold cross-validation [8, 25] to minimize the prediction error $\int_0^\tau (N_i(t) - \hat{E}N_i(t))^2 \, d\{\sum_{k=1}^n N_k(t)\}$, where $\hat{E}N_i(t) = \int_0^t Y_i(u) \exp\{\hat{\beta}(W_i)^T Z_i(u) + \hat{g}(W_i)\} \hat{\lambda}_0(u) \, du$ is the estimate of the expected failure number up to time $t$. We chose $K = 20$. Here, examination of the resulting estimated coefficient functions and their 95% confidence bands (not presented here) suggests that variable treatment and marital status are not very significant. We therefore applied the variable selection technique to the data with $\lambda = 0.02$ and bandwidth $h = 15$. The coefficient function for the treatment effect was estimated as zero at 89.5% of grid points, the coefficient function for marital status was estimated as zero at 97.9% of grid points and they were simultaneously estimated as zero at 87.5% of grid points. Thus, the variables treatment indicator and marital status were deleted. In other words, there is no significant treatment effect even when the more objective model (less restrictive model than in [29] and [18]) is used. Applying the local log partial-likelihood method (2.4) to the remaining five variables, we obtained estimated coefficient functions as in Figure 4 above. These functions depict the extent to which the gender effect and the health effect vary with age, and indicate clearly that the risk of staying at a nursing home depends on age.

## APPENDIX A: PROOFS

**A.1. Notation and conditions.** For easy reference, we collect a set of notation and conditions to be used. Let $(\Omega, \mathcal{F}, P_{(\boldsymbol{\beta}, g, \lambda)})$ be a family of complete probability spaces provided with a history $\mathbf{F} = \{\mathcal{F}_t\}_t$ for an increasing right-continuous filtration $\mathcal{F}_t \subset \mathcal{F}$. We assume that $W_i$ is $\mathcal{F}_t$-measurable, and $N_i(u)$ and $\mathbf{Z}_i(u)$ are $\mathbf{F}$-adapted. Write $\mathcal{F}_t = \sigma\{\mathbf{X}_i \leq u, \mathbf{Z}_i(u), W_i, Y_i(u), i = 1, 2, \ldots, n, 0 \leq u \leq t\}$ and $M_i(t) = N_i(t) - \int_0^t \lambda_i(u)\,du$, $i = 1, 2, \ldots, n$. Obviously, $M_i(t)$ is an $\mathcal{F}_t$ martingale.

Let $\|\cdot\|$ denote the $L_2$-norm and let $\|\cdot\|_J$ be the sup-norm of a function or a process on a set $J$. The support of the random variable $W$ is denoted by $\mathcal{W}$. For a compact subset $J_W$ of $\mathcal{W}$, we define the neighborhood set of $J_{W,\varepsilon}$ as

$$J_{W,\varepsilon} = \left\{ w \colon \inf_{w_0 \in J_W} |w - w_0| \leq \varepsilon \right\}$$

for some $\varepsilon > 0$.

To facilitate technical arguments, we will reparametrize the local partial likelihood (2.4) via the transformation $\boldsymbol{\zeta} = H(\boldsymbol{\xi} - \boldsymbol{\xi}_0)$. Hence, the logarithm of the local partial-likelihood function is

$$\tilde{\ell}_n(t, \boldsymbol{\zeta}) = \ell_n(H^{-1}\boldsymbol{\zeta} + \boldsymbol{\xi}_0, t)$$

$$= \frac{1}{n} \sum_{i=1}^{n} \int_0^t K_h(W_i - w_0)$$

$$\times [\boldsymbol{\zeta}^T \mathbf{U}_i^*(u) + \boldsymbol{\xi}_0^T \mathbf{X}_i^*(u) - \log S_{n0}(u, \boldsymbol{\zeta}, w_0)]\,dM_i(u)$$

$$+ \frac{1}{n} \sum_{i=1}^{n} \int_0^t K_h(W_i - w_0)[\boldsymbol{\zeta}^T \mathbf{U}_i^*(u) + \boldsymbol{\xi}_0^T \mathbf{X}_i^*(u) - \log S_{n0}(u, \boldsymbol{\zeta}, w_0)]$$

$$\times Y_i(u)\exp(\boldsymbol{\beta}_0(W_i)^T \mathbf{Z}_i(u) + g_0(W_i))\lambda_0(u)\,du,$$

where $\mathbf{U}_i^*(u) = \mathbf{H}^{-1}\mathbf{X}_i^*(u)$ and

$$S_{nk}(u, \boldsymbol{\zeta}, w_0) = \sum_{i=1}^n K_h(W_i - w_0)Y_i(u)\exp(\boldsymbol{\zeta}^T\mathbf{U}_i^*(u) + \boldsymbol{\xi}_0^T\mathbf{X}_i^*(u))(\mathbf{U}_i^*(u))^{\otimes k},$$

$$k = 0, 1, 2.$$

Furthermore, for each $u \in [0, \tau]$ and $k = 0, 1, 2$, we write $\tilde{\ell}_n(\boldsymbol{\zeta}) = \tilde{\ell}(\boldsymbol{\zeta}, \tau)$ and define

$$S_{nk}^*(u, \boldsymbol{\theta}, w_0) = \sum_{i=1}^n K_h(W_i - w_0)Y_i(u)\exp(\boldsymbol{\beta}^T(W_i)\mathbf{Z}_i(u) + g(W_i))(\mathbf{U}_i^*(u))^{\otimes k},$$

where $\boldsymbol{\xi}(\cdot) = (\boldsymbol{\beta}^T(\cdot), \boldsymbol{\beta}'(\cdot)^T, g(\cdot))^T$, $\boldsymbol{\theta}(\cdot) = (\boldsymbol{\beta}^T(\cdot), g(\cdot))^T$ and $w_0 \in J_W$.

Let $f(w_0)$ be the density of the random variable $W$. In addition to the notation introduced before Theorem 2, we also define, for $w_0 \in J_{W,\varepsilon}$,

$$s_0^*(u, \boldsymbol{\theta}, w_0) = f(w_0)E[\rho(u, \mathbf{Z}(u), w_0)|W = w_0],$$

$$s_1^*(u, \boldsymbol{\theta}, w_0) = f(w_0)E[\rho(u, \mathbf{Z}(u), w_0)(\mathbf{Z}^T(u), 0, 0)^T|W = w_0],$$

$$s_2^*(u, \boldsymbol{\theta}, w_0) = f(w_0)E\left[\rho(u, \mathbf{Z}(u), w_0)\exp(\boldsymbol{\beta}(w_0)^T\mathbf{Z}(u) + g(w_0))\right.$$

$$\left. \times \begin{pmatrix} \mathbf{Z}(u)\mathbf{Z}^T(u) & \mathbf{0} & \mathbf{0} \\ \mathbf{0} & \mathbf{Z}(u)\mathbf{Z}^T(u)\mu_2 & \mathbf{Z}(u)\mu_2 \\ \mathbf{0} & \mathbf{Z}^T(u)\mu_2 & \mu_2 \end{pmatrix} \middle| W = w_0\right]$$

and

$$s_k(u, \boldsymbol{\zeta}, w_0)$$

$$= f(w_0)\int E[P(u, \mathbf{Z}(u), w_0)\Psi(\boldsymbol{\zeta}, \boldsymbol{\xi}_0, \mathbf{Z}(u), y)\mathbf{R}_u(y)^{\otimes k}|W = w_0]K(y)\,dy,$$

where $k = 0, 1, 2$, $\mathbf{R}_u(y) = (\mathbf{Z}^T(u), \mathbf{Z}^T(u)y, y)^T$ and

$$\Psi(\boldsymbol{\zeta}, \boldsymbol{\xi}_0, \mathbf{Z}, y) = \exp\left(\boldsymbol{\zeta}^T\mathbf{R}_u(y) + \boldsymbol{\xi}_0^T\begin{pmatrix} \mathbf{Z} \\ \mathbf{0} \\ \mathbf{0} \end{pmatrix}\right).$$

To facilitate notation, the arguments $\boldsymbol{\theta}_0(w) = (\boldsymbol{\beta}_0^T(w), g_0(w))^T$, $\boldsymbol{\xi}_0(w)$, $\boldsymbol{\zeta}_0 = 0$ and $w_0$ are omitted in $S_{nk}^*(t, \boldsymbol{\theta}, w_0)$, $S_{nk}(t, \boldsymbol{\zeta}, w_0)$, $s_k^*(t, \boldsymbol{\theta}, w_0)$ and $s_k(t, \boldsymbol{\zeta}, w_0)$ whenever there is no ambiguity. For example,

$$S_{nk}^*(t) = S_{nk}^*(t, w_0) = S_{nk}^*(t, \boldsymbol{\theta}_0, w_0), \qquad s_k^*(t) = s_k^*(t, w_0) = s_k^*(t, \boldsymbol{\theta}_0, w_0),$$

$$S_{nk}(t) = S_{nk}(t, w_0) = S_{nk}(t, 0, w_0), \qquad s_k(t) = s_k(t, w_0) = s_k(t, 0, w_0),$$

$$S_{nk}(t, \boldsymbol{\zeta}) = S_{nk}(t, \boldsymbol{\zeta}, w_0), \qquad\qquad s_k(t, \boldsymbol{\zeta}) = s_k(t, \boldsymbol{\zeta}, w_0).$$

CONDITION A.

1. The kernel function $K \geq 0$ is a bounded, symmetric density function with compact support.
2. The functions $\boldsymbol{\beta}(\cdot)$ and $g(\cdot)$ have continuous second-order derivatives around the point $w_0$.
3. The density function $f(\cdot)$ of $W$ is continuous at the point $w_0$ and $f(w_0) > 0$.
4. The conditional probability $P(u, \mathbf{Z}(u), \cdot)$ is equicontinuous at $w_0$ and the covariate $\mathbf{Z}(u)$ is continuous.
5. We have $nh \to \infty$ and $nh^5$ is bounded.
6. We have $\int_0^\tau \lambda_0(t)\,dt < \infty$.
7. (Lindeberg condition) There exists $\delta > 0$ such that
$$(nh)^{-1/2} \sup_{t \in [0,\tau], i \in \mathcal{N}} |\mathbf{Z}_i(t)| Y_i(t) I(\boldsymbol{\beta}_0^T(w_0)\mathbf{Z}_i(t) > -\delta|\mathbf{Z}_i(t)|) \xrightarrow{P} 0,$$
where $\mathcal{N} = \{1, 2, \ldots, n\}$.
8. (Asymptotic variance) The matrix $\mathbf{a}_2 - \int_0^\tau \frac{\mathbf{a}_1(u)\mathbf{a}_1(u)^T}{\mathbf{a}_0(u)}\,d\Lambda_0(u)$ is positive definite at the point $w_0$ and the matrix $\begin{pmatrix} \mathbf{a}_2 & \mathbf{a}_1 \\ \mathbf{a}_1^T & \mathbf{a}_0 \end{pmatrix}$ is nonsingular at the point $w_0$.

Condition A will be used to derive the pointwise convergence properties of $\hat{\boldsymbol{\xi}}$ and its asymptotic normality. Conditions A.1–A.5 are similar to those in [16] and Conditions A.7–A.8 are similar to Conditions C and D of [2]. Condition A.7 seems complicated, but can be easily verified in some important cases. For example, when the covariates $\mathbf{Z}$ are bounded, the condition is always satisfied; if the covariates $\mathbf{Z}$ are bounded by a random variable that has a bounded $r$th moment for some constant $r > 2$, the condition also holds. Other cases can be found in [2]. To derive the uniformly consistent result, Condition A needs to be strengthened as follows.

CONDITION B.

1. The kernel function $K \geq 0$ is a bounded, symmetric density function with compact support.
2. The functions $\boldsymbol{\beta}_0(\cdot)$ and $g_0(\cdot)$ have continuous second-order derivatives on $J_{W,\varepsilon}$.
3. The conditional probability $P(u, \mathbf{Z}(u), w)$ is equicontinuous in the arguments $(u, w)$ on $[0, \tau] \times J_{W,\varepsilon}$.
4. The compact set $J_W \subset \mathcal{W}$ has the property $\inf_{w \in J_{W,\varepsilon}} f(w) > 0$ for some $\varepsilon > 0$.
5. The covariate process $\mathbf{Z}(u)$ has continuous sample paths in a subset $\mathcal{Z}$ of the continuous function space, and $\int_0^\tau \lambda_0(t)\,dt < \infty$ and $\|f_W\|_{J_W} < \infty$.

6. The function $s_0(t, \boldsymbol{\theta}, w_0)$ is bounded away from 0 on the product space $[0, \tau] \times \mathbf{C} \times J_{W,\varepsilon}$, that is,

$$\inf_{t \in [0,\tau]} \inf_{(\boldsymbol{\beta}^T, g) \in \mathbf{C}} \inf_{w_0 \in J_{W,\varepsilon}} s_0(t, \boldsymbol{\theta}, w_0) > 0$$

and

$$\sup_{t \in [0,\tau]} \sup_{(\boldsymbol{\beta}^T, g) \in \mathbf{C}} E|\mathbf{Z}(t)|^k \exp(\boldsymbol{\beta}^T \mathbf{Z}(t) + g) < \infty,$$

where $\mathbf{C} \subset \mathbf{R}^{p+1}$.
7. We have $nh/\log n \to \infty$ and $nh^5$ is bounded.
8. (Asymptotic variance) The matrix $\mathbf{a}_2 - \int_0^\tau \frac{\mathbf{a}_1(u)\mathbf{a}_1(u)^T}{\mathbf{a}_0(u)} \, d\Lambda_0(u)$ is positive definite for any $w_0 \in J_{W,\varepsilon}$ and the matrix $\left(\begin{smallmatrix} \mathbf{a}_2 & \mathbf{a}_1 \\ \mathbf{a}_1^T & \mathbf{a}_0 \end{smallmatrix}\right)$ is nonsingular for every $w_0 \in J_{W,\varepsilon}$.

**A.2. Proof of main results.** Let

$$C_n(t) = n^{-1} \sum_{i=1}^n Y_i(t) g(W_i, (W_i - w_0)/h, \mathbf{Z}_i(t)) K_h(W_i - w_0)$$

for a function $g(\cdot, \cdot, \cdot)$.

LEMMA A.1. *Assume that Conditions* A.1 *and* A.4 *hold. Suppose that* $g(\cdot, \cdot, \cdot)$ *is continuous in its three arguments and that* $E(g(W, u, \mathbf{Z}(t))|W = w_0)$ *is continuous at the point* $w_0$. *If* $h \to 0$ *in such a way that* $nh/\log n \to \infty$, *then*

$$\sup_{0 \le t \le \tau} |C_n(t) - C(t)| \xrightarrow{P} 0,$$

*where* $C(t) = f(w_0) \int E(Y(t)g(w_0, u, \mathbf{Z}(t))|W = w_0) K(u) \, du$.

PROOF. It is easy to show that for every $t \in [0, \tau]$,

(A.1)                          $|C_n(t) - C(t)| \xrightarrow{P} 0$.

Now we divide $[0, \tau]$ into $M$ subintervals $[t_{i-1}, t_i], i = 1, 2, \ldots, M$, with maximum length $\delta$. Then

(A.2)                          $\max_{1 < i \le M} |C_n(t_i) - C(t_i)| \xrightarrow{P} 0$.

Note that

$$\sup_{0 \le t \le \tau} |C_n(t) - C(t)|$$

(A.3)
$$\le \max_{1 \le i \le M} |C_n(t_i) - C(t_i)|$$

$$+ \max_{1 \le i \le M} \sup_{|t - t_{i-1}| < \delta} |C_n(t) - C(t) - (C_n(t_{i-1}) - C(t_{i-1}))|.$$

The first term on the right-hand side is asymptotically negligible. We now deal with the second term. Write

$$g(W, (W - w_0)/h, \mathbf{Z}) = g^+(W, (W - w_0)/h, \mathbf{Z}) - g^-(W, (W - w_0)/h, \mathbf{Z}),$$

where $g^+(\cdot, \cdot, \cdot)$ and $g^-(\cdot, \cdot, \cdot)$ are the positive part and negative part of $g(\cdot, \cdot, \cdot)$, respectively. Correspondingly, we decompose $C_n(t)$ into $C_n^+(t)$ and $C_n^-(t)$. We only need to show that

$$
\begin{aligned}
\max_{1 \le i \le M} \sup_{|t - t_{i-1}| < \delta} |C_n^+(t) - C_n^+(t_{i-1})| + \max_{1 \le i \le M} \sup_{|t - t_{i-1}| < \delta} |C^+(t) - C^+(t_{i-1})| \\
\xrightarrow{P} 0
\end{aligned}
\tag{A.4}
$$

and a similar result for $C_n^-(t)$. We now focus on (A.4). It will be shown in Appendix B that

$$\max_{1 \le i \le M} \sup_{|t - t_{i-1}| < \delta} |C_n^+(t) - C_n^+(t_{i-1})| \xrightarrow{P} 0. \tag{A.5}$$

On the other hand, we have

$$
\begin{aligned}
\max_{1 \le i \le M} & \sup_{|t - t_{i-1}| < \delta} |C^+(t) - C^+(t_{i-1})| \\
& \le \max_{1 \le i \le M} \sup_{|t - t_{i-1}| \le \delta} f(w_0) \int E\{Y(t)[g^+(w_0, u, \mathbf{Z}(t)) \\
& \qquad\qquad\qquad\qquad\qquad - g^+(w_0, u, \mathbf{Z}(t_{i-1}))]|W = w_0\} \\
& \qquad\qquad\qquad\qquad\qquad \times K(u)\, du \\
& \quad + \max_{1 \le i \le M} \sup_{|t - t_{i-1}| \le \delta} \left| \int E\{I(t_{i-1} < X < t_i)g^+(w_0, u, \mathbf{Z}(t_{i-1}))|W = w_0\} \right. \\
& \qquad\qquad\qquad\qquad\qquad\qquad \left. \times K(u)\, du \right|,
\end{aligned}
\tag{A.6}
$$

which tends to zero as $\delta \to 0$. Hence (A.4) holds. This completes the proof. $\square$

Lemma A.2. *Assume that $g(w, u, \mathbf{Z}(t))$ is equicontinuous in its arguments $w$ and $u$, and that $E(g(w_0, u, \mathbf{Z}(t))|W = w_0)$ is equicontinuous in the argument $w_0$. Under Conditions B.3 and 4, we have*

$$\sup_{0 \le t \le \tau} \sup_{w_0 \in B} |C_n(t, w_0) - C(t, w_0)| \xrightarrow{P} 0,$$

*where $B$ is a compact set that satisfies $\inf_{w \in B} f(w) > 0$.*

The proof of Lemma A.2 is similar to that of Lemma A.1 and is omitted.

LEMMA A.3. *Let $C$ and $D$ be compact sets in $R^d$ and $R^p$, and let $f(x, \boldsymbol{\theta})$ be a continuous function in $\boldsymbol{\theta} \in C$ and $x \in D$. Assume that $\boldsymbol{\theta}_0(x)$ is continuous in $x \in D$ and is the unique maximizer of $f(x, \boldsymbol{\theta})$. Let $\hat{\boldsymbol{\theta}}_n(x) \in C$ be a maximizer of $f_n(x, \boldsymbol{\theta})$. If*

$$\sup_{\boldsymbol{\theta} \in C, x \in D} |f_n(x, \boldsymbol{\theta}) - f(x, \boldsymbol{\theta})| \longrightarrow 0,$$

*then*

$$\sup_{x \in D} |\hat{\boldsymbol{\theta}}_n(x) - \boldsymbol{\theta}_0(x)| \longrightarrow 0.$$

The proof of Lemma A.3 can be found in [11].

LEMMA A.4. *Under Condition A, we have for $k = 0, 1, 2$,*

$$n^{-1} S_{nk}^*(u) = s_k^*(u) + o_p(1),$$

*uniformly for $u \in (0, \tau]$, where $s_k^*(u) = s_k^*(u, \boldsymbol{\theta}_0, w_0)$ and*

$$\sup_{u \in (0, \tau]} \|n^{-1} S_{nk}^*(u, \boldsymbol{\theta}, w_0) - s_k^*(u, \boldsymbol{\theta}, w_0)\| = o_p(1),$$

*where $\boldsymbol{\theta}$ lies in a neighborhood of $\boldsymbol{\theta}_0$ for fixed $w_0$. In addition, we have for each $\boldsymbol{\zeta}$,*

$$\sup_{u \in (0, \tau]} \|n^{-1} S_{nk}(u, \boldsymbol{\zeta}, w_0) - s_k(u, \boldsymbol{\zeta}, w_0)\| = o_p(1),$$

*where $\boldsymbol{\zeta}$ lies in a neighborhood of $\mathbf{0}$ for fixed $w_0$. Furthermore, under Condition B, we have*

$$\|n^{-1} S_{nk}^* - s_k^*\|_{\mathcal{R}} = o_p(1),$$

*where $\mathcal{R} = [0, \tau] \times \mathbf{C} \times J_{W, \varepsilon}$ and a similar result holds for $S_{nk}(u, \boldsymbol{\zeta}, w_0)$.*

The results of Lemma A.4 can be easily proved along similar lines to the arguments establishing Lemma A.1.

PROOF OF THEOREM 1. The first result of Theorem 1 follows from the first step in the proof of Theorem 2. Now we only prove the second result of Theorem 1. By an argument similar to that in the first step in the proof of Theorem 2, we easily prove from Lemma A.2 that

$$\sup_{t \in [0, \tau]} \sup_{\boldsymbol{\xi}_0 \in \mathbf{C}^*} \sup_{w_0 \in J_W} |\tilde{\ell}_n(t, \boldsymbol{\zeta}) - \tilde{\ell}_n(t, 0) - Y(t, \boldsymbol{\zeta})| \longrightarrow 0$$

in probability; here $\mathbf{C}^*$ is a convex and compact set of $R^{2p+1}$. Therefore, it follows from Lemma A.3 that $\sup_{w_0 \in J_W} |\hat{\boldsymbol{\zeta}}| \to 0$ in probability. The proof is complete. $\quad\square$

PROOF OF THEOREM 2. We first prove that $\sqrt{nh}\mathbf{H}(\hat{\boldsymbol{\xi}}(w_0) - \boldsymbol{\xi}_0(w_0))$ is asymptotically normal with mean $h^2\mathbf{e}_p\boldsymbol{\xi}_0''(w_0)\mu_2/2$ and covariance $\boldsymbol{\Sigma}(\tau, w_0)$. Now we divide the proof of the asymptotic normality of $\sqrt{nh}\mathbf{H}(\hat{\boldsymbol{\xi}}(w_0) - \boldsymbol{\xi}_0(w_0))$ into three steps. The first step is to show that $\mathbf{H}(\hat{\boldsymbol{\xi}}(w_0) - \boldsymbol{\xi}_0(w_0)) \to 0$ in probability. The second step is to establish the asymptotic normality of the first derivative of the local partial likelihood. The third step is to demonstrate that the Hessian matrix of the local partial-likelihood function converges to a positive definite one. Theorem 2 will then be proved by combining the results in these three steps.

(a) We first show that $\hat{\boldsymbol{\zeta}} \to 0$ in probability, where $\hat{\boldsymbol{\zeta}} = \mathbf{H}(\hat{\boldsymbol{\xi}} - \boldsymbol{\xi}_0)$. It is easy to show that

$$
\begin{aligned}
(A.7) \quad & \tilde{\ell}_n(t, \boldsymbol{\zeta}) - \tilde{\ell}_n(t, \mathbf{0}) \\
& = \frac{1}{n}\sum_{i=1}^n \int_0^t K_h(W_i - w_0)\left[\boldsymbol{\zeta}^T\mathbf{U}_i^*(u) - \log\frac{S_{n0}(u,\boldsymbol{\zeta})}{S_{n0}(u,0)}\right]dM_i(u) \\
& \quad + \frac{1}{n}\int_0^t S_{n1}^*(u)^T\boldsymbol{\zeta}\lambda_0(u)\,du - \frac{1}{n}\int_0^t \log\frac{S_{n0}(u,\boldsymbol{\zeta})}{S_{n0}(u,0)}S_{n0}^*(u)\lambda_0(u)\,du \\
& := X_n(t, \boldsymbol{\zeta}) + Y_n(t, \boldsymbol{\zeta}).
\end{aligned}
$$

By Lemma A.1 we obtain that

$$
\begin{aligned}
Y_n(t, \boldsymbol{\zeta}) & = \int_0^t (s_1^*(u))^T\boldsymbol{\zeta}\lambda_0(u)\,du - \int_0^t \log\frac{s_0(u,\boldsymbol{\zeta})}{s_0(u,0)}s_0^*(u)\lambda_0(u)\,du + o_p(1) \\
& := Y(t, \boldsymbol{\zeta}) + o_p(1).
\end{aligned}
$$

In Appendix B, we will show that $Y(t, \boldsymbol{\zeta})$ is a strictly concave function in $\boldsymbol{\zeta}$ and has maximum value at $\boldsymbol{\zeta} = 0$. The process $X_n(t, \boldsymbol{\zeta})$ is a local square integrable martingale with the square variation process

$$
\begin{aligned}
D_n(t) & = \langle X_n(\cdot, \boldsymbol{\zeta}), X_n(\cdot, \boldsymbol{\zeta})\rangle(t) \\
& = \frac{1}{n^2}\sum_{i=1}^n \int_0^t K_h^2(W_i - w_0)\left[\boldsymbol{\zeta}^T\mathbf{U}_i^*(u) - \log\left(\frac{S_{n0}(u,\boldsymbol{\zeta})}{S_{n0}(u,0)}\right)\right]^{\otimes 2} \\
& \quad \times Y_i(u)\exp(\boldsymbol{\beta}_0(W_i)^T\mathbf{Z}_i(u) + g_0(W_i))\lambda_0(u)\,du.
\end{aligned}
$$

It follows from Lemma A.1 that

$$
EX_n^2(t, \boldsymbol{\zeta}) = ED_n(t) = O((nh)^{-1}) \longrightarrow 0, \qquad 0 < t \le \tau.
$$

Hence, we have that

$$\tilde{\ell}_n(t, \boldsymbol{\zeta}) - \tilde{\ell}_n(t, \mathbf{0}) = Y(t, \boldsymbol{\zeta}) + O_p((nh)^{-1/2}).$$

Obviously, $\tilde{\ell}_n(t, \boldsymbol{\zeta}) - \tilde{\ell}_n(t, \mathbf{0})$ is strictly concave in $\boldsymbol{\zeta}$ with the maximizer $\hat{\boldsymbol{\zeta}}$. By the concavity lemma it follows that $\hat{\boldsymbol{\zeta}} \to 0$, the maximizer of $Y(t, \boldsymbol{\zeta})$ in probability.

(b) We now show that $\sqrt{nh}(\tilde{\ell}'_n(\tau, \mathbf{0}) - \mathbf{B}_n(\tau, w_0))$ is asymptotically normal with mean zero and covariance $\boldsymbol{\Sigma}(\tau, w_0)$, where the definitions of $\mathbf{B}_n(\tau, w_0)$ and $\boldsymbol{\Sigma}(\tau, w_0)$ can be found below.

Observe that

$$\tilde{\ell}'_n(0) = \tilde{\ell}'_n(\tau, 0) = \frac{1}{n} \sum_{i=1}^{n} \int_0^\tau K_h(W_i - w_0) \left[ \mathbf{U}_i^*(u) - \frac{S_{n1}(u, w_0)}{S_{n0}(u, w_0)} \right] dM_i(u)$$

$$+ \frac{1}{n} \sum_{i=1}^{n} \int_0^\tau K_h(W_i - w_0) \left[ \mathbf{U}_i^*(u) - \frac{S_{n1}(u, w_0)}{S_{n0}(u, w_0)} \right]$$

$$\times \exp(\boldsymbol{\beta}_0(W_i)^T \mathbf{Z}_i(u) + g_0(W_i)) Y_i(u) \lambda_0(u) \, du.$$

Let us denote the above two terms, respectively, by $I_1(\tau, 0)$ and $I_2(\tau, 0)$. We first deal with $I_2(\tau, 0)$. By Taylor expansion we have

$$\begin{aligned}
&\exp(\boldsymbol{\beta}_0(W_i)^T \mathbf{Z}_i(u) + g_0(W_i)) - \exp(\boldsymbol{\xi}_0^T \mathbf{X}_i^* + g_0(w_0)) \\
\text{(A.8)} \quad &= \tfrac{1}{2} \exp(\boldsymbol{\xi}_0^T \mathbf{X}_i^* + g_0(w_0))[\boldsymbol{\beta}_0''(w_0)^T \mathbf{Z}_i(u) + g_0''(w_0)] \\
&\times (W_i - w_0)^2 (1 + O_p(h)).
\end{aligned}$$

Note that

$$\begin{aligned}
I_2(\tau, 0) &= \frac{1}{n} \sum_{i=1}^{n} \int_0^\tau K_h(W_i - w_0) \left( \mathbf{U}_i^*(u) - \frac{S_{n1}(u)}{S_{n0}(u)} \right) \\
&\times [\exp(\boldsymbol{\beta}_0(W_i)^T \mathbf{Z}_i(u) + g_0(W_i)) - \exp(\boldsymbol{\xi}_0^T \mathbf{X}_i^* + g_0(w_0))] \\
&\times Y_i(u) \lambda_0(u) \, du.
\end{aligned}$$

Then it follows from Lemmas A.1 and A.4 that

$$\begin{aligned}
I_2(\tau, 0) &= \frac{1}{2n} \sum_{i=1}^{n} \int_0^\tau K_h(W_i - w_0) \left[ \mathbf{U}_i^*(u) - \frac{s_1^*(u)}{s_0^*(u)} \right] \\
&\times Y_i(u) \exp(\boldsymbol{\xi}_0^T \mathbf{X}_i^* + g_0(w_0))[\boldsymbol{\beta}_0''(w_0)^T \mathbf{Z}_i(u) + g_0''(w_0)] \\
&\times (W_i - w_0)^2 \lambda_0(u) \, du (1 + o_p(h)) \\
&= \frac{1}{2} h^2 f(w_0) \int_0^\tau E \left\{ \left[ \begin{pmatrix} \mathbf{Z}(u)\mu_2 \\ \mathbf{Z}(u)\mu_3 \\ \mu_3 \end{pmatrix} - \frac{s_1^*(u)\mu_2}{s_0^*(u)} \right] \rho(u, \mathbf{Z}(u), w_0) \right.
\end{aligned}$$

$$\times \; [\beta_0''(w_0)^T \mathbf{Z}(u) + g_0''(w_0)]|W = w_0 \Bigg\}$$

$$\times \; \lambda_0(u) \, du (1 + O_p(h)),$$

where $s_k^*(u) = s_k^*(u, \boldsymbol{\theta}_0, w_0)$ for $k = 0, 1, 2$. Since $K(\cdot)$ is a symmetric function, which implies $\mu_3 = 0$, simple algebra shows that

$$
\begin{aligned}
I_2(\tau, 0) &= \tfrac{1}{2} h^2 \mu_2 f(w_0) \\
&\quad \times \int_0^\tau E \left[ \begin{pmatrix} \mathbf{Z}(u) - \mathbf{a}_1(u)/\mathbf{a}_0(u) \\ 0 \\ 0 \end{pmatrix} \right. \\
&\quad \left. \times \rho(u, \mathbf{Z}(u), w_0)(\mathbf{Z}^T, 0, 1) \begin{pmatrix} \boldsymbol{\beta}_0''(w_0) \\ 0 \\ g''(w_0) \end{pmatrix} \right] d\Lambda_0(u) \\
&\quad \times (1 + O_p(h)) \\
&= \tfrac{1}{2} h^2 \mu_2 \mathbf{e}_p \boldsymbol{\Gamma}^{-1} \boldsymbol{\beta}_0''(w_0)(1 + O_p(h)).
\end{aligned}
\tag{A.9}
$$

Let us denote the term in (A.9) by $\mathbf{B}_n(\tau, w_0)$.

We now derive the asymptotic normality of the term $I_1(\tau, 0)$. Let $I_1^*(t) = \sqrt{nh} I_1(t, 0)$. Then

$$
\begin{aligned}
\langle I_1^*, I_1^* \rangle(t) &= \frac{h}{n} \sum_{i=1}^n \int_0^t K_h^2(W_i - w_0) \left[ \mathbf{U}_i^*(u) - \frac{S_{n1}(u)}{S_{n0}(u)} \right]^{\otimes 2} \\
&\quad \times Y_i(u) \exp(\boldsymbol{\beta}_0(W_i)^T \mathbf{Z}_i(u) + g_0(W_i)) \lambda_0(u) \, du.
\end{aligned}
$$

By Lemma A.1 and using Conditions A.1 and A.8, it can be shown that

$$
\begin{aligned}
\boldsymbol{\Pi}(\tau, w_0) &= \lim_{n \to \infty} E \langle I_1^*, I_1^* \rangle(\tau) \\
&= f(w_0) \\
&\quad \times \int_0^\tau E \left[ \begin{pmatrix} (\mathbf{Z}(u) - \mathbf{a}_1(u)/\mathbf{a}_0(u))^{\otimes 2} \nu_0 & 0 & 0 \\ 0 & \mathbf{Z}(u)\mathbf{Z}^T(u)\nu_2 & \mathbf{Z}(u)\nu_2 \\ 0 & \mathbf{Z}^T(u)\nu_2 & \nu_2 \end{pmatrix} \right. \\
&\qquad\qquad \left. \times \rho(u, \mathbf{Z}(u), w_0) \middle| W = w_0 \right] d\Lambda_0(u) \\
&= \begin{pmatrix} \boldsymbol{\Gamma}^{-1}\nu_0 & 0 & 0 \\ 0 & \mathbf{a}_2 \nu_2 & \mathbf{a}_1 \nu_2 \\ 0 & \mathbf{a}_1^T \nu_2 & \mathbf{a}_0 \nu_2 \end{pmatrix}.
\end{aligned}
\tag{A.10}
$$

By Condition A.7 and a proof similar to that of Anderson and Gill [2], it is easy to prove that the Lindeberg condition for the process $I_1^*(t)$ holds. By the martingale central limit theorem, we derive that $I_1^*(t)$ is asymptotically normal with mean zero and covariance $\mathbf{\Pi}(t, w_0)$. Hence

$$(A.11) \qquad \sqrt{nh}(\tilde{\ell}'_n(\mathbf{0}) - \mathbf{B}_n(\tau, w_0)) \longrightarrow N(0, \mathbf{\Pi}(\tau, w_0)).$$

(c) We will show that the second derivative of the logarithm of the local partial-likelihood function converges to a finite constant matrix. Since $\hat{\boldsymbol{\zeta}} \to 0$ in probability, by the mean-value theorem we have that

$$(A.12) \qquad \tilde{\ell}''_n(\hat{\boldsymbol{\zeta}}) = \tilde{\ell}''_n(\mathbf{0}) + o_p(1).$$

Since $s_k^*(u) = s_k(u)\exp(g_0(w_0)), k = 0, 1, 2$, from Lemma A.4, we can obtain

$$\tilde{\ell}''_n(\mathbf{0}) = \frac{1}{n}\int_0^\tau \sum_{i=1}^n K_h(W_i - w_0)\frac{s_2^*(u)s_0^*(u) - s_1^*(u)(s_1^*(u))^T}{(s_0^*(u))^2}\,dN_i(u) + o_p(1).$$

Write $F_w(u) = P(X \le u, \Delta = 1|W = w)$ and its corresponding empirical conditional measure,

$$\tilde{F}_w(u) = \frac{1}{n}\sum_{i=1}^n K_h(W_i - w_0)I(X_i \le u, \Delta_i = 1).$$

By kernel smoothing techniques, we easily prove that

$$(A.13) \qquad \begin{aligned} \tilde{\ell}''_n(0) &= -\int_0^\tau \frac{s_2^*(u)s_0^*(u) - s_1^*(u)(s_1^*(u))^T}{(s_0^*(u))^2}\,d\tilde{F}_w(u) + o_p(1) \\ &= -\mathbf{A}(\tau, w_0) + o_p(1), \end{aligned}$$

where

$$\mathbf{A} = \int_0^\tau \frac{s_2^*(u)s_0^*(u) - s_1^*(u)(s_1^*(u))^T}{(s_0^*(u))^2}\,dF_w(u).$$

It is easy to show that $\mathbf{A}(\tau, w_0)$ is positive definite.

(d) Combining the results in steps (a), (b) and (c), we can establish the asymptotic normality of $\sqrt{nh}\mathbf{H}(\hat{\boldsymbol{\xi}}(w_0) - \boldsymbol{\xi}(w_0))$. In fact, since $\hat{\boldsymbol{\zeta}}$ maximizes $\tilde{\ell}_n(\boldsymbol{\zeta})$, by Taylor expansion around 0, we have

$$-\tilde{\ell}'_n(\mathbf{0}) = \tilde{\ell}'_n(\hat{\boldsymbol{\zeta}}) - \tilde{\ell}'_n(\mathbf{0}) = (\tilde{\ell}''_n(\hat{\boldsymbol{\zeta}}^*))^T\hat{\boldsymbol{\zeta}},$$

where $\hat{\boldsymbol{\zeta}}^*$ lies between 0 and $\hat{\boldsymbol{\zeta}}$. Hence $\hat{\boldsymbol{\zeta}}^* \to 0$ in probability. It follows from (A.13) that

$$\begin{aligned} &\hat{\boldsymbol{\zeta}} - \mathbf{A}(\tau, w_0)^{-1}\mathbf{B}_n(\tau, w_0) \\ &= -(\tilde{\ell}''_n(\hat{\boldsymbol{\zeta}}^*))^{-1}(\tilde{\ell}'_n(\mathbf{0}) - \mathbf{B}_n(\tau, w_0)) + o_p(1). \end{aligned}$$

Combining (A.11) with (A.13), by Slutsky's theorem we obtain that

$$\sqrt{nh}(\hat{\boldsymbol{\zeta}} - \mathbf{A}(\tau, w_0)^{-1} \mathbf{B}_n(\tau, w_0))$$
$$\longrightarrow N(\mathbf{0}, \mathbf{A}^{-1}(\tau, w_0) \boldsymbol{\Pi}(\tau, w_0)(\mathbf{A}^{-1}(\tau, w_0))^T).$$

Now we simplify the matrix $\mathbf{A}(\tau, w_0)$. Obviously, by a simple calculation we have

$$s_2^*(u) = f(w_0) E \left[ \begin{pmatrix} \mathbf{Z}(u)\mathbf{Z}^T(u) & \mathbf{0} & \mathbf{0} \\ \mathbf{0} & \mathbf{Z}(u)\mathbf{Z}^T(u)\mu_2 & \mathbf{Z}(u)\mu_2 \\ \mathbf{0} & \mathbf{Z}^T(u)\mu_2 & \mu_2 \end{pmatrix} \right.$$

$$(A.14) \qquad\qquad\qquad\qquad \left. \times \, \rho(u, \mathbf{Z}(u), w_0) | W = w_0 \right]$$

$$= \begin{pmatrix} \mathbf{a}_2(u) & \mathbf{0} & \mathbf{0} \\ \mathbf{0} & \mathbf{a}_2(u)\mu_2 & \mathbf{a}_1(u)\mu_2 \\ \mathbf{0} & \mathbf{a}_1^T(u)\mu_2 & \mathbf{a}_0(u)\mu_2 \end{pmatrix}.$$

Similarly, we obtain that

$$(A.15) \qquad (s_1^*(u))^{\otimes 2} = \begin{pmatrix} \mathbf{a}_1(u)\mathbf{a}_1^T(u) & \mathbf{0} & \mathbf{0} \\ \mathbf{0} & \mathbf{0} & \mathbf{0} \\ \mathbf{0} & \mathbf{0} & 0 \end{pmatrix}.$$

Note that $s_0^*(u) = \mathbf{a}_0(u)$. Hence it follows from (A.14) and (A.15) that

$$(A.16) \qquad \mathbf{A}(\tau, w_0) = \begin{pmatrix} \boldsymbol{\Gamma}^{-1} & \mathbf{0} & \mathbf{0} \\ \mathbf{0} & \mathbf{a}_2\mu_2 & \mathbf{a}_1\mu_2 \\ \mathbf{0} & \mathbf{a}_1^T\mu_2 & \mathbf{a}_0\mu_2 \end{pmatrix}.$$

Hence, the asymptotic bias of the estimator $\hat{\boldsymbol{\zeta}}(w_0)$ is

$$b(\tau, w_0) = \mathbf{A}^{-1}(\tau, w_0) \mathbf{B}_n(\tau, w_0)$$
$$= h^2 \mathbf{e}_p \boldsymbol{\xi}_0''(w_0) \mu_2 / 2$$

and the asymptotic covariance is

$$\boldsymbol{\Sigma}(\tau, w_0) = \mathbf{A}^{-1}(\tau, w_0) \boldsymbol{\Pi}(\tau, w_0)(\mathbf{A}^{-1}(\tau, w_0))^T$$
$$= \begin{pmatrix} \boldsymbol{\Gamma}\nu_0 & \mathbf{0} \\ \mathbf{0}^T & \begin{pmatrix} \mathbf{a}_2 & \mathbf{a}_1 \\ \mathbf{a}_1 & \mathbf{a}_0 \end{pmatrix}^{-1} \mu_2^{-2}\nu_2 \end{pmatrix}$$
$$= \begin{pmatrix} \boldsymbol{\Gamma} & \mathbf{0} \\ \mathbf{0}^T & \mathbf{Q}\mu_2^{-2}\nu_2 \end{pmatrix}.$$

This completes the proof. □

PROOF OF THEOREM 3. We have shown from (A.13) that

$$\tilde{\ell}_n''(\hat{\boldsymbol{\zeta}}^*) = -\mathbf{A}(\tau, w_0) + o_p(1) \tag{A.17}$$

for any $\hat{\boldsymbol{\zeta}}^*$ between zero and $\hat{\boldsymbol{\zeta}} = \mathbf{H}(\hat{\boldsymbol{\xi}} - \boldsymbol{\xi}_0)$. By Theorem 2, $\hat{\boldsymbol{\zeta}} = O_p(h^2 + (nh)^{-1/2})$. Thus, for any $\hat{\boldsymbol{\zeta}}^* = O_p(h^2 + (nh)^{-1/2})$, (A.17) holds.

By Taylor expansion of $\tilde{\ell}_n'(\hat{\boldsymbol{\zeta}}_0)$ at $\boldsymbol{\zeta}_0 = 0$, we have

$$\tilde{\ell}_n'(\hat{\boldsymbol{\zeta}}_0) = \tilde{\ell}_n'(\boldsymbol{\zeta}_0) + \tilde{\ell}_n''(\hat{\boldsymbol{\zeta}}^*)(\hat{\boldsymbol{\zeta}}_0 - \boldsymbol{\zeta}_0), \tag{A.18}$$

where $\hat{\boldsymbol{\zeta}}_0 = \mathbf{H}(\hat{\boldsymbol{\xi}}_0 - \boldsymbol{\xi}_0)$ and $\hat{\boldsymbol{\zeta}}^* = \mathbf{H}(\hat{\boldsymbol{\xi}}_0^* - \boldsymbol{\xi}_0)$, in which $\hat{\boldsymbol{\xi}}_0^*$ lies between $\boldsymbol{\xi}_0$ and $\hat{\boldsymbol{\xi}}_0$.

By definition of the one-step estimator and (A.18), we have that

$$\hat{\boldsymbol{\zeta}}_{\mathrm{os}} - \boldsymbol{\zeta}_0 = (\hat{\boldsymbol{\zeta}}_0 - \boldsymbol{\zeta}_0) - (\tilde{\ell}_n''(\hat{\boldsymbol{\zeta}}_0))^{-1} \tilde{\ell}_n'(\hat{\boldsymbol{\zeta}}_0).$$

Using (A.17), we have

$$\begin{aligned}
\hat{\boldsymbol{\zeta}}_{\mathrm{os}} - \boldsymbol{\zeta}_0 &= (I - (\tilde{\ell}_n''(\hat{\boldsymbol{\zeta}}_0))^{-1} \tilde{\ell}_n''(\hat{\boldsymbol{\zeta}}^*))(\hat{\boldsymbol{\zeta}}_0 - \boldsymbol{\zeta}_0) - (\tilde{\ell}_n''(\hat{\boldsymbol{\zeta}}_0))^{-1} \tilde{\ell}_n'(\boldsymbol{\zeta}_0) \\
&= -(\tilde{\ell}_n''(\hat{\boldsymbol{\zeta}}_0))^{-1} \tilde{\ell}_n'(\boldsymbol{\zeta}_0) + o_p(\hat{\boldsymbol{\zeta}}_0 - \boldsymbol{\zeta}_0) \\
&= -(\tilde{\ell}_n''(\hat{\boldsymbol{\zeta}}_0))^{-1} [\tilde{\ell}_n'(\boldsymbol{\zeta}_0) - \mathbf{B}_n(\tau, w_0)] - \tilde{\ell}_n''(\hat{\boldsymbol{\zeta}}_0)^{-1} \mathbf{B}_n(\tau, w_0) \\
&\quad + o_p((nh)^{-1/2} + h^2).
\end{aligned}$$

It follows from (A.11) and (A.13) that $\hat{\boldsymbol{\zeta}}_{\mathrm{os}}$ has the same asymptotic distribution as the maximum local partial-likelihood estimator. This yields Theorem 3. □

PROOF OF THEOREM 4. By the same argument as that of Lemma A.1, we have

$$\sup_{t \in [0, \tau]} \sup_{\|\boldsymbol{\theta} - \boldsymbol{\theta}_0\| \le \|\hat{\boldsymbol{\theta}} - \boldsymbol{\theta}_0\|} n^{-1} |\Delta_n(t, \boldsymbol{\theta}) - \Delta_n(t, \boldsymbol{\theta}_0)| \longrightarrow 0 \tag{A.19}$$

in probability, where

$$\Delta_n(t, \boldsymbol{\theta}) = \sum_{i=1}^n I(W_i \in J_W) Y_i(t) \exp\{\boldsymbol{\beta}^T(W_i)\mathbf{Z}_i(t) + g(W_i)\},$$

where $\boldsymbol{\theta} = (\boldsymbol{\beta}^T(\cdot), g(\cdot))^T$.

By definition of $\hat{\Lambda}_0(t)$, we have

$$\begin{aligned}
\hat{\Lambda}_0(t) - \Lambda_0(t) &= \int_0^t \left\{ \frac{1}{\Delta_n(\hat{\boldsymbol{\theta}})} - \frac{1}{\Delta_n(\boldsymbol{\theta}_0)} \right\} d\bar{N}_n + \int_0^t \left\{ \frac{d\bar{N}_n}{\Delta_n(\boldsymbol{\theta}_0)} - d\Lambda_0 \right\} \\
&= -\int_0^t \frac{\Delta_n(\hat{\boldsymbol{\theta}}) - \Delta_n(\boldsymbol{\theta}_0)}{\Delta_n(\hat{\boldsymbol{\theta}})} d\Lambda_0 - \int_0^t \frac{\Delta_n(\hat{\boldsymbol{\theta}}) - \Delta_n(\boldsymbol{\theta}_0)}{\Delta_n(\hat{\boldsymbol{\theta}})\Delta_n(\boldsymbol{\theta}_0)} d\bar{M}_n \\
&\quad + \int_0^t \frac{1}{\Delta_n(\boldsymbol{\theta}_0)} d\bar{M}_n,
\end{aligned}$$

where $\bar{N}_n = \sum_{i=1}^n N_i$ and $\bar{M}_n = \sum_{i=1}^n M_i$. From (A.19) it is easy to see that the first term converges to zero in probability uniformly on $(0, \tau]$ as $n \to \infty$. The last two terms of the above expression are square integrable local martingales with variation processes

$$\int_0^t \frac{(\Delta_n(\hat{\boldsymbol{\theta}}) - \Delta_n(\boldsymbol{\theta}_0))^2}{(\Delta_n(\hat{\boldsymbol{\theta}}))^2 \Delta_n(\boldsymbol{\theta}_0)} \, d\Lambda_0 \quad \text{and} \quad \int_0^t \frac{1}{\Delta_n(\boldsymbol{\theta}_0)} \, d\Lambda_0,$$

respectively. Since $\Delta_n(\boldsymbol{\theta}_0) = O_p(n)$, the above variance processes converge to zero in probability uniformly on $(0, \tau]$ as $n \to \infty$. The terms converge to zero in probability uniformly on $(0, \tau]$ by an argument similar to that of Andersen and Gill [2] via the Lenglart inequality. Therefore

$$\hat{\Lambda}_0(t) \longrightarrow \Lambda_0(t)$$

uniformly on $(0, \tau]$. Thus, we can prove by the standard argument of kernel estimation that

$$\hat{\lambda}_0(t) \longrightarrow \lambda_0(t)$$

uniformly on $(0, \tau]$. $\quad \square$

PROOF OF THEOREM 5. From the proof of Theorem 2, we easily show that this theorem holds. $\quad \square$

PROOF OF THEOREM 6. Using the same proof as in Theorem 2, we can get

$$\ell_n'(\boldsymbol{\xi}_0) = O_p((nh)^{-1/2} + h^2).$$

Let $\alpha_n = (nh)^{-1/2} + h^2 + a_n$. Following the same lines as the proof of Theorem 1 of [17], the result follows. $\quad \square$

LEMMA A.5. *Suppose that the conditions of Theorem 6 hold. Then with probability tending to 1, for any given $\boldsymbol{\xi}_1$ satisfying $\|\boldsymbol{\xi}_1 - \boldsymbol{\xi}_{10}\| = O_p((nh)^{-1/2} + h^2)$ and any constant $C$, we have*

$$Q((\boldsymbol{\xi}_1^T, 0)^T) = \max_{\|\boldsymbol{\xi}_2\| \leq C[(nh)^{-1/2} + h^2]} Q((\boldsymbol{\xi}_1^T, \boldsymbol{\xi}_2^T)^T).$$

PROOF. From an argument similar to that in step (b) in the proof of Theorem 2, it is easy to show that

$$\ell_n'(\boldsymbol{\xi}_0) = O_p((nh)^{-1/2} + h^2),$$

and by an argument similar to that in step (c) of the proof of Theorem 2, we have

$$\ell_n''(\boldsymbol{\xi}_0) = O_p(1).$$

The result follows from the the proof of Lemma 1 of [17].  □

PROOF OF THEOREM 7. It follows from Lemma A.5 and Theorem 6 that the first result of Theorem 7 holds. Now we prove the second result of Theorem 7. It can be easily shown that there exists a $\hat{\boldsymbol{\xi}}_1$ as in Theorem 6 that is a local maximizer of $Q(\boldsymbol{\xi}_1^T, \mathbf{0})^T$, and that satisfies the likelihood equations

$$\frac{\partial Q(\boldsymbol{\xi})}{\partial \boldsymbol{\xi}_1}\bigg|_{\boldsymbol{\xi}=(\hat{\boldsymbol{\xi}}_1, \mathbf{0})} = 0.$$

Using the Taylor expansion of $(\partial Q(\boldsymbol{\xi}))/\partial \boldsymbol{\xi}_1$ at point $\boldsymbol{\xi}_0$ and noting that $\hat{\boldsymbol{\xi}}_1$ is a consistent estimator from Theorem 6, we have

$$
\begin{aligned}
(\text{A.20}) \quad & \frac{\partial \ell_n(\boldsymbol{\xi}_0)}{\partial \boldsymbol{\xi}_1} + \left(\frac{\partial^2 l_n(\boldsymbol{\xi}_0)}{\partial \boldsymbol{\xi}_1 \partial \boldsymbol{\xi}_1^T} + o_p(1)\right)(\hat{\boldsymbol{\xi}}_1 - \boldsymbol{\xi}_{10}) \\
& - b - (\boldsymbol{\Sigma}_1 + o_p(1))(\hat{\boldsymbol{\xi}}_1 - \boldsymbol{\xi}_{10}) = 0.
\end{aligned}
$$

From the proof of Theorem 2, it is easy to show that

$$\sqrt{nh}\left(\mathbf{H}^{-1}\frac{\partial \ell_n(\boldsymbol{\xi}_0)}{\partial \boldsymbol{\xi}} - \frac{1}{2}h^2\mu_2\mathbf{e}_p\boldsymbol{\Gamma}^{-1}\boldsymbol{\beta}''(w_0)(1 + o_p(1))\right) \longrightarrow N(0, \boldsymbol{\Pi}(\tau, w_0))$$

and

$$\mathbf{H}^{-1}\frac{\partial^2 \ell_n(\boldsymbol{\xi}_0)}{\partial \boldsymbol{\xi}\,\partial \boldsymbol{\xi}^T}\mathbf{H}^{-1} \longrightarrow -\mathbf{A}(\tau, w_0).$$

Thus, we have

$$
\begin{aligned}
(\text{A.21}) \quad & \sqrt{nh}\left(\mathbf{H}_1^{-1}\frac{\partial \ell_n(\boldsymbol{\xi}_0)}{\partial \boldsymbol{\xi}_1} - \frac{1}{2}h^2\mu_2\begin{pmatrix}\boldsymbol{\Gamma}_{-1}\boldsymbol{\beta}_0''(w_0) \\ 0\end{pmatrix}(1 + o_p(1))\right) \\
& \longrightarrow N(0, \boldsymbol{\Pi}_1(\tau, w_0))
\end{aligned}
$$

and

$$(\text{A.22}) \qquad \mathbf{H}_1^{-1}\frac{\partial^2 l_n(\boldsymbol{\xi}_0)}{\partial \boldsymbol{\xi}_1 \partial \boldsymbol{\xi}_1^T}\mathbf{H}_1^{-1} \longrightarrow -\mathbf{A}_1(\tau, w_0).$$

By some simple calculations, we easily show that the second result of Theorem 7 follows from (A.20), (A.21) and (A.22).  □

## APPENDIX B

*Concavity and maxima of* $Y(t, \beta)$. Here we prove that $Y(\tau, \boldsymbol{\zeta})$ defined by (A.7) is concave with respect to $\boldsymbol{\zeta}$. Differentiating the function $Y(\tau, \boldsymbol{\zeta})$ with respect to $\boldsymbol{\zeta}$, we have

$$\frac{\partial Y(\tau,\boldsymbol{\zeta})}{\partial \boldsymbol{\zeta}} = \int_0^t s_1^*(u)\lambda_0(u)\,du - \int_0^t \frac{s_1(u,\boldsymbol{\zeta})}{s_0(u,\boldsymbol{\zeta})}s_0^*(u)\lambda_0(u)\,du,$$

$$\frac{\partial^2 Y(\tau,\boldsymbol{\zeta})}{\partial \boldsymbol{\zeta}^2} = -\int_0^t \frac{s_2(u,\boldsymbol{\zeta})s_0(u,\boldsymbol{\zeta}) - (s_1(u,\boldsymbol{\zeta}))^{\otimes 2}}{(s_0(u,\boldsymbol{\zeta}))^2}s_0^*(u)\lambda_0(u)\,du.$$

By the integral transform and the fact that $\mathbf{aa^T} + \mathbf{bb^T} \geq \mathbf{2ab^T}$ for any vectors $\mathbf{a}$ and $\mathbf{b}$, we can show that

$$\frac{\partial^2 Y(\tau,0)}{\partial \boldsymbol{\zeta}^2} < 0.$$

Again by $s_k^*(u,0) = s_k(u,\boldsymbol{\theta}_0)\exp(g_0(w_0)), k = 0,1,$ we have

$$\frac{\partial Y(\tau,0)}{\partial \boldsymbol{\zeta}} = 0.$$

Hence $\boldsymbol{\zeta} = 0$ is the maximizer $Y(\tau,\boldsymbol{\zeta})$.

PROOF OF (A.5). It is easy to show that

$$\max_{1 \leq i \leq M}\sup_{|t-t_{i-1}|<\delta}|C_n^+(t) - C_n^+(t_{i-1})| \leq J_1 + J_2,$$

where

$$J_1 = \max_{1 \leq i \leq M}\sup_{|t-t_{i-1}|\leq\delta}\left| n^{-1}\sum_{j=1}^n Y_j(t)g^+(W_j,(W_j-w_0)/h,\mathbf{Z}_j(t))K_h(W_j-w_0) \right.$$
$$- n^{-1}\sum_{j=1}^n Y_j(t)g^+(W_j,(W_j-w_0)/h,\mathbf{Z}_j(t_{i-1}))$$
$$\left. \times K_h(W_j-w_0) \right|$$

and

$$J_2 = \max_{1 \leq i \leq M}\sup_{|t-t_{i-1}|\leq\delta}\left| n^{-1}\sum_{j=1}^n Y_j(t)g^+(W_j,(W_j-w_0)/h,\mathbf{Z}_j(t_{i-1})) \right.$$
$$\times K_h(W_j-w_0)$$
$$- n^{-1}\sum_{i=1}^n Y_j(t_{i-1})g^+(W_j,(W_j-w_0)/h,\mathbf{Z}_j(t_{i-1}))$$
$$\left. \times K_h(W_j-w_0) \right|.$$

Note that $\mathbf{Z}_j(t)$ $(j = 1, 2, \ldots, n)$ is continuous on $[0, \tau]$. Thus we easily obtain that

$$
\begin{aligned}
J_1 \leq \max_{1 \leq j \leq n} \sup_{t \in [0,\tau]} \sup_{|t - t_{i-1}| \leq \delta} &|g^+(W_j, (W_j - w_0)/h, \mathbf{Z}_j(t) \\
&- g^+(W_j, (W_j - w_0)/h, \mathbf{Z}_j(t_{i-1}))| \\
&\times \sup_{t \in [0,\tau]} n^{-1} \sum_{j=1}^{n} Y_j(t) K_h(W_j - w_0),
\end{aligned}
$$

which tends to zero in probability. Since $Y_i(t)$ is a decreasing function of $t$, we have, for any $\varepsilon > 0$,

$$
\begin{aligned}
P(J_2 > \varepsilon) \leq M P \Bigg( n^{-1} \Bigg| \sum_{j=1}^{n} I(t_{i-1} < X_j < t_i) \\
\times g^+(W_j, (W_j - w_0)/h, \mathbf{Z}_j) K_h(W_j - w_0) \Bigg| > \varepsilon \Bigg).
\end{aligned}
$$

It is easy to show that

$$
n^{-1} \sum_{j=1}^{n} I(t_{i-1} < X_j < t_i) g^+(W_j, (W_j - w_0)/h, \mathbf{Z}_j(t_{i-1})) K_h(W_j - w_0)
$$

$$
\xrightarrow{P} f(w_0) \int E\{I(t_{i-1} < X < t_i) g^+(w_0, u, \mathbf{Z}_j(t_{i-1}) | W = w_0)\} K(u)\, du.
$$

On the other hand,

$$
\begin{aligned}
E(I(t_{i-1} < X < t_i) &g^+(w_0, u, \mathbf{Z}(t_{i-1})) | W = w_0) \\
&\leq E^{1/2} \{ I(t_{i-1} < X < t_i) | W = w_0) \} \\
&\quad \times E^{1/2} \{ g^{+2}(w_0, u, \mathbf{Z}(t_{i-1})) | W = w_0 \} \\
&= |P(X < t_{i-1} | W = w_0) - P(X < t_i | W = w_0)|^{1/2} \\
&\quad \times E^{1/2}(g^{+2}(w_0, u, \mathbf{Z}(t_{i-1})) | W = w_0) \\
&< \varepsilon
\end{aligned}
$$

as $|t_i - t_{i-1}| < \delta$. Hence

$$
\begin{aligned}
P \Bigg( n^{-1} \Bigg| \sum_{j=1}^{n} I(t_{i-1} < X_j < t_i) \\
\times g^+(W_j, (W_j - w_0)/h, \mathbf{Z}_j(t_{i-1})) K_h(W_j - w_0) \Bigg| > \varepsilon \Bigg)
\end{aligned}
$$

$$\leq P\Bigg(\Bigg|\sum_{j=1}^{n} n^{-1} I(t_{i-1} < X_j < t_i)$$

$$\times g^+(W_j, (W_j - w_0)/h, \mathbf{Z}_j(t_{i-1})) K_h(W_j - w_0)$$

$$- f(w_0) \int E(I(t_{i-1} < X < t_i)$$

$$\times g^+(w_0, u, \mathbf{Z}(t_{i-1}))|W = w_0) K(u)\, du\Bigg| > \varepsilon/2\Bigg)$$

$$+ P\Bigg(f(w_0) \int |E(I(t_{i-1} < X < t_i)$$

$$\times g^+(w_0, u, \mathbf{Z}(t_{i-1}))|W = w_0) K(u)\, du| > \varepsilon/2\Bigg)$$

$$< \eta.$$

Hence for any $\eta > 0$ and $\varepsilon > 0$ there exists $N_0$ such that for $n > N_0$ we have

$$(B.1) \qquad\qquad P(J_1 + J_2 > \varepsilon) < 2\eta.$$

Therefore, we obtain that

$$(B.2) \qquad P\Bigg(\max_{1 \leq i \leq M} \sup_{|t - t_{i-1}| < \delta} |C_n^+(t) - C_n^+(t_{i-1})| > \varepsilon\Bigg) < 2\eta.$$

This completes the proof of (A.5). $\qquad\square$

**Acknowledgment.** The authors are grateful to the reviewers for constructive comments that have significantly improved the presentation of the paper.

J. FAN
DEPARTMENT OF STATISTICS
CHINESE UNIVERSITY OF HONG KONG
SHATIN
HONG KONG
AND
DEPARTMENT OF OPERATIONS RESEARCH
    AND FINANCIAL ENGINEERING
PRINCETON UNIVERSITY
PRINCETON, NEW JERSEY 08540
USA
E-MAIL: jqfan@princeton.edu

H. LIN
SCHOOL OF MATHEMATICS
SICHUAN UNIVERSITY
CHENGDU, SICHUAN 610064
PEOPLE'S REPUBLIC OF CHINA

Y. ZHOU
INSTITUTE OF APPLIED MATHEMATICS
ACADEMY OF MATHEMATICS AND SYSTEM SCIENCE
CHINESE ACADEMY OF SCIENCE
ZHONG GUANCUN, BEIJING 100080
PEOPLE'S REPUBLIC OF CHINA
E-MAIL: yzhou@amss.ac.cn